\documentclass[12pt, reqno]{amsart}
\usepackage[latin1]{inputenc}
\usepackage{amsmath}
\usepackage{amsfonts}
\usepackage{amssymb}

\newcommand{\barr}{\bar}

\newcommand{\Aut}{\text{Aut}_{\mathbb{C}}(\mathcal{M},0)}
\newcommand{\C}{\mathbb{C}}
\newcommand{\M}{\mathcal{M}}
\newcommand{\BB}{\mathcal{B}}
\newcommand{\CC}{\mathcal{C}}
\newcommand{\HH}{\mathcal{H}}
\newcommand{\R}{\mathbb{R}}

\newcommand{\X}{\chi}
\newcommand{\T}{\tau}
\newcommand{\f}{\tilde{f}}
\newcommand{\g}{\tilde{g}}

\theoremstyle{plain}
\newtheorem{aaaa}{Theorem}[section]
\newtheorem{cor1.2}[aaaa]{Corollary}
\theoremstyle{definition}
\newtheorem{remark2.2.4}[aaaa]{Remark}
\theoremstyle{plain}
\newtheorem{cor1.3}[aaaa]{Corollary}

\theoremstyle{plain}
\newtheorem{theorem1.4}[aaaa]{Theorem}

\theoremstyle{definition}
\newtheorem{example2}[aaaa]{Example}
\theoremstyle{plain}
\newtheorem{algebraic}[aaaa]{Theorem}

\theoremstyle{definition}
\newtheorem{segreClaim}{Proposition}[section]

\newtheorem{theorem1.1}{Theorem}[section]
\newtheorem{conditionD}[theorem1.1]{Definition}
\newtheorem{example3}[theorem1.1]{Example}

\theoremstyle{plain}
\newtheorem{finitejet}[theorem1.1]{Theorem}
\theoremstyle{definition}
\newtheorem{remark3}[theorem1.1]{Remark}
\theoremstyle{plain}
\newtheorem{lemma2.2}[theorem1.1]{Lemma}
\newtheorem{lemma2.3}[theorem1.1]{Lemma}
\newtheorem{lemma2.4}[theorem1.1]{Lemma}
\theoremstyle{definition}
\newtheorem{defn2.5}[theorem1.1]{Definition}
\theoremstyle{plain}
\newtheorem{lemma2.6}[theorem1.1]{Lemma}
\newtheorem{lemma2.7}[theorem1.1]{Lemma}
\newtheorem{lemma2.8}[theorem1.1]{Lemma}
\newtheorem{lemma2.9}[theorem1.1]{Lemma}
\newtheorem{lemma2.10}[theorem1.1]{Lemma}
\newtheorem{theorem2.11}[theorem1.1]{Theorem}
\newtheorem{algebraicII}[theorem1.1]{Theorem}
\newtheorem{lemma3.1}[theorem1.1]{Lemma}

\theoremstyle{plain}
\newtheorem{lemma4.4}[theorem1.1]{Lemma}
\theoremstyle{definition}
\newtheorem{segSubObs}[theorem1.1]{Observation}
\newtheorem{nextObs}[theorem1.1]{Observation}
\newtheorem{finitetype}{Theorem}[section]
\newtheorem{example31}[finitetype]{Example}
\newtheorem{example32}[finitetype]{Example}

\newtheorem{ex4.1}[finitetype]{Example}
\newtheorem{ex4.2}[finitetype]{Example}
\newtheorem{ex4.3}[finitetype]{Example}
\newtheorem{ex4.4}[finitetype]{Example}

\numberwithin{equation}{section}
\begin {document}

\def\1#1{\ov{#1}}
\def\2#1{\widetilde{#1}}
\def\3#1{\mathcal{#1}}
\def\4#1{\widehat{#1}}

\title[Parametrization of Holomorphic Segre preserving Maps] {Parametrization of Holomorphic Segre preserving Maps}
\author[R. B. Angle]{R. Blair Angle}
\address{ Department of Mathematics, University of California
at San Diego, La Jolla, CA 92093-0112, USA}
\email{angle@metsci.com }

\maketitle

\begin{abstract}
In this paper, we explore holomorphic Segre preserving maps.
First, we investigate holomorphic Segre preserving maps
  sending the complexification $\mathcal{M}$ of a generic real analytic
 submanifold $M \subseteq \C^N$
of finite type at some point $p$ into the complexification $\mathcal{M}'$ of a  generic real analytic
 submanifold $M' \subseteq \C^{N'}$, finitely nondegenerate at some point $p'$.
We prove that for a fixed $M$ and $M'$, the  germs
at $(p,\bar{p})$ of Segre submersive holomorphic Segre preserving maps
 sending $(\M,(p,\bar{p}))$ into $(\M',(p', \bar{p}'))$ can be parametrized by their $r$-jets at $(p,\bar{p})$, for some fixed $r$ depending only on $M$ and $M'$.
(If, in addition, $M$ and $M'$ are both real algebraic, then we prove that any such map must  be holomorphic algebraic.)
From this  parametrization, it follows that the
set of germs of holomorphic Segre preserving automorphisms
$\mathcal{H}$ of the complexification $\mathcal{M}$ of a real
analytic submanifold finitely nondegenerate and of finite type at
some point $p$, and such that $\mathcal{H}$ fixes $(p,\bar{p})$, is
an algebraic complex Lie group.
We then explore the relationship between this automorphism group and the
group of automorphisms of $M$ at $p$.
%We find parametrizations for these
%automorphism groups for several submanifolds, some of
%which are finitely degenerate at $p$.
\end{abstract}

\section{Introduction}
Let $M \subseteq \C^{N}$ be a  real analytic submanifold of codimension $d$, with $p \in M$, given locally near $p$ by the real
analytic defining function $\rho(Z,\bar{Z})$.  The \emph{complexification} $\M$ of $M$ is a holomorphic submanifold of $\C^{2N}$ given locally for $(Z,\zeta) \in  \C^{N} \times \C^{N}$ near $(p,\bar{p})$ by
$\M = \{(Z,\zeta) : \rho(Z,\zeta)=0 \}$. Now assume $M$ is generic (see Section 2), and let $M' \subseteq \C^{N'}$ be a generic real analytic submanifold of codimension $d'$,
with $p' \in M'$, and let $\M'$ denote its complexification.
%Let $Z=(z,w) \in \C^n \times \C^d$ and $\zeta=(\X,\T) \in \C^n \times \C^d$, and
Consider a holomorphic  map $\HH: (\C^{2N},(p,\bar{p})) \rightarrow
(\C^{2N'},(p',\bar{p}'))$ defined on a neighborhood of $(p,\bar{p})$  of the form
\begin{equation}  \label{hspm}
\HH(Z,\zeta)=\big( H(Z),\widetilde{H}(\zeta)\big),
\end{equation}
where $H, \widetilde{H} : \C^N \rightarrow \C^{N'}$.
 Assume further that $\HH(\M) \subseteq \M'$. These maps will be the chief object
of study in this paper.  We will call such a map a \emph{holomorphic Segre preserving map} (HSPM) as it preserves Segre varieties in a
sense which will be made  precise in Section \ref{section:background}.
 Utilizing the notation $\overline{\varphi}(z):= \overline{\varphi(\bar{z})}$, we observe that if  $\widetilde{H}=\overline{H}$,
 then $H$ is a holomorphic map defined near
 $p$ sending $(M,p)$ into $(M',p')$. Such maps have been extensively studied. However, HSPMs
are  relatively new and unstudied objects (for related recent work, see  \cite{Angle}, \cite{A07b}, and \cite{Z07}).
Under certain
restrictions, the collection of HSPMs sending $\M$ into $\M'$ is,
%in a sense to be made more precise
in a manner to be described in more detail
in subsequent sections, ``bigger" than the
collection of holomorphic mappings sending $M$ into $M'$. We shall see several examples of this in
Section \ref{section:Examples}.

For $p_0 \in \C^m$, let $\mathcal{T}_{p_0}(\C^m)$ denote the holomorphic tangent space of $\C^m$ at $p_0$.
Let $\mathcal{T}_{(p,\bar{p})}^0 \M  \subseteq \mathcal{T}_{(p,\bar{p})}(\C^{2N})$ denote the set of all vectors of the form $\sum_{j=1}^N a_j \frac{\partial}{\partial Z_j}
+ \sum_{j=1}^N b_j \frac{\partial}{\partial \zeta_j}$ such that
$\sum_{j=1}^N a_j \frac{\partial}{\partial Z_j}$
and ${\sum_{j=1}^N b_j \frac{\partial}{\partial \overline{Z}_j}}$ are tangent to $M$ at $p$.
A vector of the form $\sum_{j=1}^N a_j \frac{\partial}{\partial Z_j}$ tangent to $M$ at $p$ is known as a \emph{holomorphic tangent vector},
and a vector of the form $\sum_{j=1}^N b_j \frac{\partial}{\partial \overline{Z}_j}$
tangent to $M$ at $p$ is known as an \emph{antiholomorphic tangent vector}.
For any HSPM $\HH$ sending $(\M,(p,\bar{p}))$ into $(\M',(p',\bar{p}'))$,
 $\displaystyle{\mathcal{D}_{(p,\bar{p})}\HH \left(\mathcal{T}_{(p,\bar{p})}^0 \M \right) \subseteq \mathcal{T}_{(p',\bar{p}')}^0 \M'}$, where $\mathcal{D}_{(p,\bar{p})}\HH : \mathcal{T}_{(p,\bar{p})}(\C^{2N}) \rightarrow \mathcal{T}_{(p',\bar{p}')} (\C^{2N'})$ is defined by
 $\mathcal{D}_{(p,\bar{p})}\HH(\mathcal{L})(\varphi)=\mathcal{L}(\varphi \circ \HH)$ for any holomorphic function $\varphi:(\C^{2N'},(p',\bar{p}')) \rightarrow \C$.
We say that $\HH$ is \emph{Segre submersive}
at $(p,\bar{p})$ if
$$\mathcal{D}_{(p,\bar{p})}\HH \left(\mathcal{T}_{(p,\bar{p})}^0 \M \right) = \mathcal{T}_{(p',\bar{p}')}^0 \M'.$$
This definition is independent of choice of coordinates for $M$ and $M'$.

Given $M$ and $M'$ satisfying certain geometric conditions, our main result, Theorem \ref{thm:Theorem1.1},
states that the  germs
at $(p,\bar{p})$ of HSPMs, Segre submersive at $(p,\bar{p})$, sending $(\M,(p,\bar{p}))$ into $(\M',(p', \bar{p}'))$ can be parametrized by their $r$-jets, for some fixed $r$ depending only on $M$ and $M'$.
This result  was motivated by, and is a
generalization of, results due to Baouendi, Ebenfelt, and Rothschild \cite{BER99b} and
Baouendi, Rothschild, and Zaitsev \cite{BRZ01}. We also mention a recent
paper of Lamel and Mir \cite{LM07} for related results.
Before stating Theorem \ref{thm:Theorem1.1}, we present some more notation.
Let $J^K(\mathbb{C}^N,\mathbb{C}^{N'})_{(p,p')}$ denote
the set of $K$-jets  at $p$ of germs of holomorphic maps from
$(\mathbb{C}^N,p)$ into $(\mathbb{C}^{N'},p')$.
(In this paper, we assume that $J^K(\mathbb{C}^N,\mathbb{C}^{N'})_{(p,p')}$ includes only derivatives of \emph{positive} order.)
Let $j_p^K$ represent the corresponding \emph{K-jet map} defined
on the set of germs at $p$ of holomorphic mappings  given by
$$j_p^K \phi = \Bigg( \frac{\partial^{|\alpha|} \phi}{\partial Z^\alpha} (p) \Bigg)_{1 \leq |\alpha| \leq K}.$$
%Without loss of generality, in order to simplify
%notation and the statement of results, we will make three assumptions for the remainder of this paper. Unless
%otherwise stated, assume that:
%\renewcommand{\labelenumi}{(\roman{enumi})}
%\begin{enumerate}
%\item $M$ and $M'$ are always understood to be generic and real analytic.
%\item $p=p'=0$.
%\item Any HSPM $\HH$ maps 0 to 0 and is given in the form (\ref{hspm}).
%\end{enumerate}

%For the remainder of this paper, unless otherwise stated, we will assume $M$ and $M'$ are generic and real analytic, $p \in M$ and $p' \in M'$, %and any HSPM sends $(\M,(p,\bar{p}))$ into $(\M',(p',\bar{p}'))$.

\begin{aaaa} \label{thm:Theorem1.1}
Let $M \subseteq \mathbb{C}^{N}$ be real analytic, generic, and of finite type at $p$. Let $M'
\subseteq \mathbb{C}^{N'}$ be real analytic, generic, and  finitely nondegenerate at $p'$. Then
there exist positive integers $K$ and $r$, depending only on $M$ and $M'$, and $\C^{N'}$-valued
holomorphic functions $\Phi_1, \ldots, \Phi_r$ defined on an open subset of
$\mathbb{C}^{N} \times J^K(\mathbb{C}^{N},
\mathbb{C}^{N'})_{(p,p')} \times J^K(\mathbb{C}^{N},
\mathbb{C}^{N'})_{(\bar{p},\bar{p}')} $ of the form
\begin{equation}      \label{phi1}
\Phi_l(Z,\Lambda,\Gamma)=\sum_{\gamma} \frac{
P^l_\gamma(\Lambda,\Gamma)  }{ Q^l_1(\Lambda)^{s^l_\gamma}
Q^l_2(\Gamma)^{t^l_\gamma}   }  \left(Z - p \right)^{\gamma} ,
\end{equation}
where $s_\gamma^l$ and $t_\gamma^l$ are nonnegative integers,  $P_\gamma^l$ are
$\mathbb{C}^{N'}$-valued polynomials, and $Q_1^l$ and $Q_2^l$ are
$\C$-valued polynomials with real coefficients,  such that the following holds. Let
$\mathcal{H}(Z,\zeta)=\big(H(Z),\widetilde{H}(\zeta)\big)$ be  an HSPM
sending $(\M,(p,\bar{p}))$ into $(\M',(p', \bar{p}'))$
such that $\HH$ is Segre submersive at $(p,\bar{p})$. Then there exists $1 \leq l \leq r$ such that
%satisfying $Q_1^l(j_p^{K}(H)) \neq 0$,
%$Q_1^l(j_{\bar{p}}^{K}(\widetilde{H})) \neq 0$,  $Q_2^l(j_p^{K}(H)) \neq 0$, and  $Q_2^l(j_{\bar{p}}^{K}(\widetilde{H})) \neq 0$ such that
%for any such $l$,
\begin{equation}
H(Z) = \Phi_l \Big(Z, j_p^{K}(H), j_{\bar{p}}^K(\widetilde{H})   \Big) ,
\end{equation}
\begin{equation}  \label{csfd}
\widetilde{H}(\zeta) = \overline{\Phi_l} \Big(\zeta, j_{\bar{p}}^{K}(\widetilde{H}),
j_p^K(H)   \Big) ,
\end{equation}
for $(Z,\zeta)$ sufficiently close to $(p,\bar{p})$.
Furthermore, for any $(\Lambda_0, \Gamma_0) \in
J^K(\mathbb{C}^{N}, \mathbb{C}^{N'})_{(p,p')} \times
J^K(\mathbb{C}^{N}, \mathbb{C}^{N'})_{(\bar{p},\bar{p}')}$ such that
$Q_1^l(\Lambda_0) \neq 0$ and $Q_2^l(\Gamma_0) \neq 0$, $\Phi_l$ is
holomorphic in a neighborhood of $(p,\Lambda_0,\Gamma_0)$.
\end{aaaa}

The appearance of $\overline{\Phi_l}$ in (\ref{csfd}) is interesting and will be instrumental in the proof of Corollary \ref{thm:cor1.3}.
The reader is referred to Section \ref{section:background} for precise definitions of
\emph{finite type} and \emph{finite nondegeneracy}.

Define
$$\text{Aut}(M,p) := \{H:(\C^N,p) \rightarrow (\C^N,p) \, | \,   H(M) \subseteq M ,$$
\begin{equation*} \label{eq:AutReal}
 H \text{ is a germ at } p \text{ of a holomorphic map}, H \text{ is invertible at } p \} ,
\end{equation*}
and
$$\text{Aut}_{\mathbb{C}}(\mathcal{M},(p,\bar{p})) := \{\mathcal{H}:(\C^{2N},(p,\bar{p})) \rightarrow (\C^{2N},(p,\bar{p})) \, | \,  \HH(\M) \subseteq \M ,$$
\begin{equation*} \label{eq:AutSegre}
 \mathcal{H} \text{ is a germ at } (p,\bar{p}) \text{ of an HSPM } , \mathcal{H} \text{ is invertible at } (p,\bar{p}) \} .
\end{equation*}
We call $\text{Aut}(M,p)$ the \emph{group of automorphisms} of $M$ at $p$,
and we call $\text{Aut}_{\mathbb{C}}(\mathcal{M},(p,\bar{p}))$ the
\emph{group of holomorphic Segre preserving automorphisms} of $\M$ at $(p,\bar{p})$.
Let  $J_p^K(\mathbb{C}^{N}) := J^K(\mathbb{C}^{N},\mathbb{C}^{N})_{(p,p)}$ be a simplification of notation,
define $G_p^K(\mathbb{C}^{N})$ to be the set of all elements of $J_p^K(\mathbb{C}^{N})$
which correspond to invertible mappings at $p$, and define a jet map $\eta_{(p,\bar{p})}^K$ on the set of germs at $(p,\bar{p})$ of HSPMs such that for
 $\HH=(H,\widetilde{H})$, $\eta_{(p,\bar{p})}^K(\HH) := (j_p^K H, j_{\bar{p}}^K \widetilde{H})$.
Theorem \ref{thm:Theorem1.1} then leads to the following corollary.

\begin{cor1.2} \label{thm:cor1.2}
Let $M \subseteq \C^N$ be  of finite type at $p$ and finitely
nondegenerate at $p$. Then there exists an integer $K$ depending only
on $M$ such that  $\eta_{(p,\bar{p})}^K$ restricted to $\text{Aut}_{\mathbb{C}}(\mathcal{M},(p,\bar{p}))$ is a homeomorphism
onto  a closed,  holomorphic algebraic submanifold (Lie group) of
$G_p^K(\mathbb{C}^{N}) \times G_{\bar{p}}^K(\mathbb{C}^{N})$.
\end{cor1.2}

\begin{remark2.2.4} \label{remark2.2.4}
 We observe that a  consequence of Corollary \ref{thm:cor1.2} is that $j_p^K$ restricted to
 $\text{Aut}(M,p)$ is a
homeomorphism onto a closed, real algebraic submanifold
(Lie group) of $G_p^K(\mathbb{C}^{N})$.
 This fact
 has already been proven in previous work. In the case that $M$ is a hypersurface,
 it was shown by Baouendi, Ebenfelt, and Rothschild in \cite{BER97} that $j_p^K \big(\text{Aut}(M,p) \big)$
  is a closed,  real analytic
  submanifold (Lie group) of $G_p^K(\mathbb{C}^{N})$. However, it was not shown that it is also real algebraic.
  This fact was later proven for submanifolds of any codimension by Baouendi, Ebenfelt, and Rothschild in \cite{BER99b}.
  \end{remark2.2.4}

 As $j_p^K \big(\text{Aut}(M,p) \big)$ is a real algebraic
submanifold, it is natural to consider its complexification as a
holomorphic submanifold of $G_p^K(\mathbb{C}^{N}) \times
G_{\bar{p}}^K(\mathbb{C}^{N})$. We will denote this complexification
$\mathbb{C} \big\{ j_p^K \big(\text{Aut}(M,p) \big) \big\}$. As
$j_p^K \big(\text{Aut}(M,p) \big)$ is real algebraic, it has global
defining functions, and thus so does its complexification.
Similarly, $\eta_{(p,\bar{p})}^K \big(
\text{Aut}_{\mathbb{C}}(\mathcal{M},(p,\bar{p})) \big)$ has global defining functions. So a
natural question to consider is the relationship between $\mathbb{C}
\big\{ j_p^K \big(\text{Aut}(M,p) \big) \big\}$ and  $\eta_{(p,\bar{p})}^K \big(
\text{Aut}_{\mathbb{C}}(\mathcal{M},(p,\bar{p})) \big)$. The following
corollary says that the former is always contained in the latter, and
 they are necessarily of the same dimension. Does equality hold? As it turns out, sometimes there is equality, and sometimes there is not.
 In Section \ref{section:Examples} we will give
examples demonstrating both.

\begin{cor1.3} \label{thm:cor1.3}
Let $M$ and $K$ be as in Corollary \ref{thm:cor1.2}.  Let $\BB \subseteq
G_p^{K}(\mathbb{C}^{N}) \times G_{\bar{p}}^{K}(\mathbb{C}^{N})$ denote
the connected component of $\C \big\{
j_p^K\big(\text{Aut}(M,p)\big)\big\}$ which contains $(Id,Id')$,
where $Id$ (resp., $Id'$) is the point in $G_p^{K}(\mathbb{C}^{N})$ (resp., $G_{\bar{p}}^{K}(\mathbb{C}^{N})$) corresponding to the identity map on $\C^{N}$.  Let
$\CC \subseteq G_p^{K}(\mathbb{C}^{N}) \times
G_{\bar{p}}^{K}(\mathbb{C}^{N})$ denote the connected component of
$\eta_{(p,\bar{p})}^K \big(
\text{Aut}_{\mathbb{C}}(\mathcal{M},(p,\bar{p})) \big)$ which contains $(Id,Id')$. Then:
 
\begin{enumerate}
  \item $ \C \big\{ j_p^K\big(\text{Aut}(M,p)\big)\big\} \subseteq
  \eta_{(p,\bar{p})}^K \big(
\text{Aut}_{\mathbb{C}}(\mathcal{M},(p,\bar{p})) \big)$
  \item $\BB = \CC$
  \item $\eta_{(p,\bar{p})}^K \big(
\text{Aut}_{\mathbb{C}}(\mathcal{M},(p,\bar{p})) \big)
   \text{ and  } \C \big\{ j_p^K\big(\text{Aut}(M,p)\big)\big\}$  are made up of finitely many disjoint
cosets of  $\BB$.
\end{enumerate}
\end{cor1.3}

One of the strengths of Theorem \ref{thm:Theorem1.1} lies in the fact that the form
of $\Phi_l$ leads to Corollaries \ref{thm:cor1.2} and \ref{thm:cor1.3}. These functions, however,
depend upon the jets of both $H$ \emph{and} $\widetilde{H}$. In Theorem
\ref{thm:Theorem1.4}, we see that it is in fact possible, though, to find  functions
which express $\mathcal{H}$
entirely in terms of the $L$-jets of $H$ (or of $\widetilde{H}$) for some $L$. In particular, once we know $H$, we also know
$\widetilde{H}$,
and vice-versa.

\begin{theorem1.4} \label{thm:Theorem1.4}
Let $M$ and $M'$ be as in Theorem \ref{thm:Theorem1.1}. Then there exist  positive
integers $r$ and $L$, depending only on $M$ and $M'$, and $\C^{2N'}$-valued holomorphic functions
$\Phi^1_1, \ldots, \Phi^1_{r}$  defined on an open subset of $\mathbb{C}^{2N}
\times J^L(\mathbb{C}^{N},\mathbb{C}^{N'})_{(p,p')}$
and $\Phi^2_1, \ldots, \Phi^2_{r}$ defined on an open subset of $\mathbb{C}^{2N}
\times J^L(\mathbb{C}^{N},\mathbb{C}^{N'})_{(\bar{p},\bar{p}')}$
such that the following holds. Let
$\mathcal{H}(Z,\zeta)=\big(H(Z),\widetilde{H}(\zeta) \big)$ be an HSPM
sending $(\M,(p,\bar{p}))$ into $(\M',(p', \bar{p}'))$
such that $\HH$ is Segre submersive at $(p,\bar{p})$.  Then there exist $1 \leq l_1, l_2 \leq r$ such that
%for which $\Phi^1_{l_1}(Z,\zeta,j_p^L (H))$  and $\Phi^2_{l_1}(Z,\zeta,j_{\bar{p}}^L (\widetilde{H}))$ are defined for $(Z,\zeta)$ in an open %neighborhood of $(p,\bar{p})$
\begin{equation}
\mathcal{H}(Z,\zeta)= \Phi^1_{l_1} \big(Z,\zeta,j_p^L (H)\big) ,
\end{equation}
\begin{equation}
\mathcal{H}(Z,\zeta)= \Phi^2_{l_2} \big(Z,\zeta,j_{\bar{p}}^L (\widetilde{H})\big) ,
\end{equation}
for $(Z,\zeta)$ sufficiently close to $(p,\bar{p})$.
%Furthermore, for any $(\Lambda_0, \Gamma_0) \in
%J^K(\mathbb{C}^{N}, \mathbb{C}^{N'})_{(p,p')} \times
%J^K(\mathbb{C}^{N}, \mathbb{C}^{N'})_{(\bar{p},\bar{p}')}$ such that
%$Q_1(\Lambda_0) \neq 0$ and $Q_2(\Gamma_0) \neq 0$, $\Phi$ is
%holomorphic in a neighborhood of $(p,\Lambda_0,\Gamma_0)$.
\end{theorem1.4}

 Note that Theorem \ref{thm:Theorem1.4} does not necessarily hold if $M'$ is finitely degenerate at $p'$, as the following example demonstrates.

\begin{example2}
 Let $M=M' \subseteq \mathbb{C}^2$
be given by $M=\{\text{Im } w=|z|^4 \}$ and its complexification by $\mathcal{M}=\{w-\tau=2iz^2\chi^2\}$,
where $(z,w)$ and $(\X,\T)$ are coordinates on $\C^2$.
We note that $M$ is of finite type but finitely \emph{degenerate} at 0.
Let $H(z,w)=(z,w)$. We can find two \emph{distinct}
maps $\widetilde{H}_1(\chi,\tau)$ and $\widetilde{H}_2(\chi,\tau)$ such that $\mathcal{H}_1=(H,\widetilde{H}_1)$ and
$\mathcal{H}_2=(H,\widetilde{H}_2)$
both satisfy the hypotheses of Theorem \ref{thm:Theorem1.4}. Indeed, let $\widetilde{H}_1(\chi,\tau)=(\chi,\tau)$ and
let $\widetilde{H}_2(\chi,\tau)=(-\chi,\tau)$.
\end{example2}

Finally, we present a result on algebraicity. Recall that a real analytic (resp., holomorphic) mapping is said to be real analytic
 (resp., holomorphic) algebraic
if all of its components are real analytic (resp., holomorphic) algebraic, and a real analytic (resp., holomorphic) submanifold is said to be real (resp., holomorphic) algebraic
if it can be given by real analytic (resp., holomorphic) algebraic defining functions.

\begin{algebraic} \label{thm:algebraic}
Let $M$ and $M'$ be as in Theorem \ref{thm:Theorem1.1}, and assume that
$M$ and $M'$ are real algebraic. Then any HSPM
sending $(\M,(p,\bar{p}))$ into $(\M',(p', \bar{p}'))$
which is Segre submersive at $(p,\bar{p})$ is holomorphic algebraic.
\end{algebraic}

  The layout of this paper is as follows.
In Section \ref{section:2.2}, we present some additional background material. Section \ref{section:2.3} contains the reformulations
and proofs of three of the main results as given in Section 1, while Section \ref{section:2.4} is dedicated to proving the main
results of Section 1. Section \ref{section:Examples} consists of several examples of HSPMs and automorphism groups.
In particular, examples demonstrating
both  equality and non-equality  of $\mathbb{C}
\big\{ j_p^K \big(\text{Aut}(M,p) \big) \big\}$ and  $\eta_{(p,\bar{p})}^K \big(
\text{Aut}_{\mathbb{C}}(\mathcal{M},(p,\bar{p})) \big)$ are provided.
(We refer the reader to \cite{Angle} for additional examples.)
%Finally, in Section \ref{section:2.7} we present some open questions.

\section{Additional Background}  \label{section:background}   \label{section:2.2}
Let $M \subseteq \C^{N}$ be a real analytic  submanifold of codimension $d$. Recall that this means
 that given any $p \in M$, there exists a real analytic function $\rho=(\rho_1,\ldots,\rho_d):(\C^{N}, p) \rightarrow \R^d$,
 satisfying $d\rho_1 \wedge \ldots \wedge d\rho_d \neq 0$ at $p$,
 such that $M$ is given locally near $p$ by the vanishing of $\rho$.  If, in addition, $\rho$, known as the
 \emph{defining function} of $M$, satisfies the stronger condition $\partial \rho_1 \wedge \ldots \wedge \partial \rho_d \neq 0$ at
 $p$, then we say that $M$ is \emph{generic}.
If $M$ is generic, it can be shown (see, for example, \cite{BER99a}) that there exists a holomorphic change of coordinates $Z=(z,w) \in  \C^{N-d} \times \C^d$, vanishing at $p$, and an open neighborhood $\Omega$ of 0 such that in these coordinates $M$ is locally given by $\{(z,w) \in \Omega : w=Q(z,\bar{z},\bar{w}) \}$, where $Q(z,\chi,\tau)$ is a $\C^d$-valued holomorphic function defined near 0 in $\C^{N-d} \times \C^{N-d} \times \C^d$
and satisfying
%\begin{equation*}
$Q(0,\chi,\tau) \equiv Q(z,0,\tau) \equiv \tau.$
%\end{equation*}
Such coordinates are called \emph{normal coordinates}.

%Define $$\displaystyle{ \mathcal{D}_M := \Bigg\{ \sum_{j=1}^N \phi_j(Z,\zeta)\frac{\partial}{\partial Z_j}
%+ \sum_{j=1}^N \psi_j(Z,\zeta)\frac{\partial}{\partial \zeta_j} \, \Bigg| \,
%\sum_{j=1}^N \phi_j(Z,\bar{Z})\frac{\partial}{\partial Z_j}
%\text{ and } }$$ $$\displaystyle{\sum_{j=1}^N \psi_j(Z,\bar{Z})\frac{\partial}{\partial \overline{Z}_j} \text{ are vector fields tangent to } M \text{ near } p \Bigg\} .}$$
%Let $\mathcal{D}_M$ denote the set of all  vector fields of the form
%$$\sum_{j=1}^N \phi_j(Z,\zeta)\frac{\partial}{\partial Z_j}+\sum_{j=1}^N \psi_j(Z,\zeta)\frac{\partial}{\partial \zeta_j}$$
%such that
% $\sum_{j=1}^N \phi_j(Z,\bar{Z})\frac{\partial}{\partial Z_j}$ and
%$\sum_{j=1}^N \psi_j(Z,\bar{Z})\frac{\partial}{\partial \overline{Z}_j}$ are real analytic vector fields tangent to $M$ near $p$.
%Let $\mathfrak{g}_M$ denote the complex Lie algebra generated by $\mathcal{D}_M$.
A  vector field of the form $\sum_{j=1}^{N}  a_j(Z,\bar{Z})\frac{\partial}{\partial \overline{Z}_j}$
  tangent to $M$ near $p$, where $a_j$ are smooth functions on $M$,  is called a \emph{CR vector field}.
We say that $M$ is of \emph{finite type} at $p$
   (in the sense of Kohn \cite{Kohn} and Bloom and Graham \cite{BG}) if
   the CR vector fields, their complex conjugates, and all repeated commutators of these vector fields span the
   complexified tangent space of $M$ at $p$.
   %$\text{dim}_\C \, \mathfrak{g}_M(p) = 2N - d$.
  %We say that $M$ is of \emph{finite type} at $p$
   %(in the sense of Kohn \cite{Kohn} and Bloom and Graham \cite{BG})
   %if the complex Lie algebra generated by all the CR and anti-CR vector fields tangent to $M$ has dimension $N$ at $p$.
   Letting $(\rho_j)_Z := \left(\frac{\partial \rho_j}{\partial Z_1} , \ldots , \frac{\partial \rho_j}{\partial Z_N}\right)$ and
   $L^\alpha := L_1^{\alpha_1} \cdots L_m^{\alpha_m}$, where $\alpha = (\alpha_1,\ldots,\alpha_m)$ and
   $L_1, \ldots, L_m$ is a basis for the CR vector fields of $M$ near $p$,
    we say that $M$ is \emph{finitely nondegenerate} at $p$ if there exists a nonnegative integer $K$ such that
\begin{equation}  \label{finNondeg}
\text{span} \big\{ L^\alpha (\rho_j)_Z(p):|\alpha| \leq K , 1 \leq j \leq d \big\} = \mathbb{C}^{N} .
\end{equation}
We say that $M$ is \emph{k-nondegenerate} at $p$ if
$k$ is the smallest $K$ for which (\ref{finNondeg}) holds.
%We say that $M$ is \emph{k-nondegenerate} at $p$ if $k$ is the smallest integer for which this holds.
 It is not difficult to show that if $M$ is given in normal coordinates by $w=Q(z,\bar{z},\bar{w})$ then $M$ is $k$-nondegenerate
 at 0 if and only if the matrix whose rows are
 $\big( Q_{z_j \chi^\alpha}(0,0,0)\big)_{|\alpha|\leq K}$, $1 \leq j \leq N-d$, has rank $N-d$ for $K \geq k$ and rank  less than $N-d$
 for $K < k$.
%For the remainder of this paper, we will assume all submanifolds are defined near the origin and are given in normal coordinates.

Let $M \subseteq \C^N$ be a generic real analytic submanifold  such that $p \in M$, and assume that
there exists an open neighborhood $\Omega \subseteq \C^N$ such that the complexification $\M$ of $M$ is defined on $\Omega \, \times \, ^*\Omega$, where
$^*\Omega := \{\barr{Z}:Z \in \Omega\}$.
 Given any $(Z,\zeta) \in \Omega \, \times \, ^*\Omega$, we define the \emph{Segre varieties} of $M$ as follows:
\begin{equation*}
\Sigma_Z := \{\zeta\,  \in \, ^*\Omega : \rho(Z,\zeta)= 0\} ,
\end{equation*}
\begin{equation*}
\hat{\Sigma}_\zeta := \{Z\,  \in  \Omega : \rho(Z,\zeta)= 0\} ,
\end{equation*}
where $\rho(Z,\barr{Z})$ is a defining function for $M$.
Segre varieties are named for the Italian geometer Beniamino Segre who first introduced them in 1931 (\cite{Segre}).
We note here that $\M$  is sometimes referred to as the \emph{Segre family}
associated with $M$ (see, for example, \cite{Chern}, \cite{F80}).

For $(Z',\zeta')$ coordinates on $\C^{N'} \times \C^{N'}$,
let $M' \subseteq \mathbb{C}^{N'}$ be a real analytic generic submanifold, with $p' \in M'$, and denote
its complexification  by $\M'$ and its Segre varieties
 by $\Sigma_{Z'}'$ and $\hat{\Sigma}_{\zeta '}'$. Let $\mathcal{H}:\mathbb{C}^{2N}
\rightarrow \mathbb{C}^{2N'}$  be a holomorphic map defined near
$(p,\bar{p})$ sending $(\M,(p,\bar{p}))$ into $(\M',(p',\bar{p}'))$.  Furthermore, we will assume that for any $(Z,\zeta)
\in \mathcal{M}$, there exists $(Z',\zeta ') \in \mathcal{M}'$ such
that
\begin{equation}     \label{eq:qw}
\mathcal{H}\big( \{Z\} \times \Sigma_Z \big) \subseteq \{Z'\} \times \Sigma'_{Z'} \, ,
\end{equation}
\begin{equation}
 \mathcal{H}\big( \hat{\Sigma}_\zeta \times\{ \zeta\} \big) \subseteq  \hat{\Sigma}'_{\zeta '} \times \{\zeta
 '\}  \, .
\end{equation}

  \begin{segreClaim}
 $\mathcal{H}$, when restricted to $\mathcal{M}$,  is an  HSPM of the form (\ref{hspm}).
    \end{segreClaim}
  This fact was proven for hypersurfaces in \cite{F80}, but it is true for higher codimension as well.
  For the reader's convenience, we present a proof.
  \begin{proof}
 Write
 %\begin{equation}
 $\mathcal{H}(Z,\zeta)=\big(\phi_1(Z,\zeta),\phi_2(Z,\zeta)\big)$,
 %\end{equation}
   where $\phi_1$ and $\phi_2$ are $\mathbb{C}^{N'}$-valued holomorphic functions, and write
    $\bar{p}=\left(\bar{p}_1,\bar{p}_2\right) \in \C^{N-d} \times \C^d$. As $M$ is generic, it follows from
    the implicit function theorem that
   (after a possible rearrangement of coordinates) there exists a $\C^d$-valued holomorphic function $\theta$,
   satisfying $\theta(p,\bar{p}_1)=\bar{p}_2$, such that for any $Z$
   sufficiently close to $p$, $\big(Z,\bar{p}_1,\theta(Z,\bar{p}_1)\big) \in \M$.
   For any $Z$ near $p$, define $H(Z):=\phi_1\big(Z,\bar{p}_1,\theta(Z,\bar{p}_1)\big)$.
    We claim that on $\M$, $H(Z)=\phi_1(Z,\zeta)$.  This is because (\ref{eq:qw}) implies that for any $Z_0$, $\phi_1(Z_0,\zeta)$ is
    constant for all $\zeta \in \Sigma_{Z_0}$.
     A similar argument applies to $\phi_2$.
     \end{proof}

\section{Reformulations}  \label{section:generalizations}  \label{section:2.3}

In the remainder of this paper, we will assume, unless otherwise specified, that $M \subseteq \C^{m+d}$ and $M' \subseteq \C^{n+e}$ are real analytic
generic submanifolds of codimensions $d$ and $e$, respectively.
We will further assume that
$M$ is given by $w=Q(z,\bar{z},\bar{w})$, where $Z=(z,w)$ are normal coordinates, and
$M'$ is given by $w'=Q'(z',\bar{z}',\bar{w}')$, where $Z'=(z',w')$ are normal coordinates.
Thus, the complexification $\M$ (resp., $\M'$) of $M$ (resp., $M'$) is given by
$w=Q(z,\X,\T)$ (resp., $w'=Q'(z',\X',\T')$), where $\zeta=(\X,\T) \in \C^m \times \C^d$ and $\zeta'=(\X',\T') \in \C^n \times \C^e$.
Unless otherwise specified, we will assume any HSPM $\HH$ sends $(\M,0)$ into $(\M',0)$ and is given in the form
\begin{equation}
\HH(Z,\zeta) = \big( H(Z),\widetilde{H}(\zeta)\big) = \big(f(Z),g(Z),\f(\zeta),\g(\zeta)\big) ,
\end{equation}
where $f = (f^1,\ldots,f^n)$ and $\f = (\f^1,\ldots,\f^n)$ are $\C^n$-valued holomorphic functions,
$g = (g^1,\ldots,g^e)$ and $\g = (\g^1,\ldots,\g^e)$ are $\C^e$-valued holomorphic functions,
and we  write $z=(z_1,\ldots,z_m)$, $w=(w_1,\ldots,w_d)$, $z'=(z_1',\ldots,z_n')$, and $w'=(w_1',\ldots,w_e')$
(similarly for $\X$, $\T$, $\X'$, and $\T'$).

 \subsection{Reformulation of Theorem  \ref{thm:Theorem1.1}}
 %$\,$ \vspace{.06in}

We begin with a technical definition.
\begin{conditionD}
Let $M \subseteq \mathbb{C}^{m+d}$ be of codimension $d$ and $M' \subseteq \mathbb{C}^{n+ e}$
be  of codimension $e$, and
assume $m \geq n$. Let $\mathcal{H}$ be an HSPM. Let $\mu=(\mu_1,\ldots,\mu_n)$
for some $1 \leq \mu_1< \ldots<\mu_n \leq m$ and
$\nu=(\nu_1,\ldots,\nu_n)$ for some $1 \leq \nu_1< \ldots<\nu_n \leq
m$, and assume  that $\displaystyle{ \det \left( \frac{\partial
f_k}{\partial z_{\mu_l}}(0)\right)_{1\leq k,l \leq n} \neq 0}$ and
$\displaystyle{ \det \left( \frac{\partial \tilde{f}_k}{\partial
\chi_{\nu_l}}(0)\right)_{1\leq k,l \leq n} \neq 0}$. Then we say that
the map $\mathcal{H}$ satisfies \emph{condition $D_{\mu\nu}$}.
\end{conditionD}

Let us note here that any given $\mathcal{H}$ may satisfy condition
$D_{\mu\nu}$ for several different $\mu$ and $\nu$, as the following example illustrates.
\begin{example3}
Let
$M \subseteq \mathbb{C}^4$ and $M' \subseteq \mathbb{C}^3$
 be given
by
\begin{equation}
M = \big\{ \text{Im }w= |z_1|^2+2\text{Re}(z_3\bar{z}_1- z_3\bar{z}_2) - |z_2|^2 \big\} ,
\end{equation}
\begin{equation}
M' = \big\{ \text{Im }w'=|z'_1|^2+|z'_2|^2 \big\}.
\end{equation}
Note that $M$ is of finite type at 0, and $M'$ is finitely
nondegenerate at 0. Let $\mathcal{H}$ be given by
\begin{equation}
\mathcal{H}(z,w,\X,\T)=\big( z_1+z_3,z_1-z_2,w,\chi_1-\chi_2,\chi_2+\chi_3,\tau \big).
\end{equation}
Then $\mathcal{H}$ satisfies condition $D_{\mu\nu}$ for any permisssible $\mu$ and $\nu$.
That is $\mu$ can be any one of $(1,2)$, $(1,3)$, or $(2,3)$,
 as can $\nu$.
 \end{example3}

Our main theorem, from which Theorem \ref{thm:Theorem1.1} follows, is Theorem
\ref{thm:Finitejet}. Before we present it, we introduce some notation.
Given an HSPM $\HH$,  we can write
\begin{equation}
 j_0^K H = \Big(
\big(f^j_{z_l}(0)\big)_{1 \leq l \leq m, 1 \leq j \leq n}, (j_0^K)'
H\Big),
\end{equation}
 where $(j_0^K)' H$ represents  the remaining derivatives of
$H$ at 0. Given any $\Lambda \in J_0^K(\C^{m+d},$
$\C^{n+e})_{(0,0)}$, we
will then write
\begin{equation}
\Lambda = \big( (\Lambda^{j,l})_{1 \leq l \leq m,
1 \leq j \leq n}, {\Lambda}' \big),
\end{equation}
 where $\big( j_0^K H)^{j,l}$ is
exactly $f_{z_l}^j(0)$. We define a similar decomposition for $j_0^K
\widetilde{H}$. This notation will be used several times in this paper.

\begin{finitejet} \label{thm:Finitejet}
Let $M \subseteq \mathbb{C}^{m+d}$ be of codimension $d$ and of finite type at $0$.  Let $M'
\subseteq \mathbb{C}^{n+ e}$ be of codimension $e$ and $k$-nondegenerate at $0$. Then
there exists a positive integer $K$ depending only on $M$ and $M'$ such that
for each $\alpha=(\alpha_1,\ldots,\alpha_n)$ with $1\leq \alpha_1 <
\ldots < \alpha_n \leq m$ and each $\beta=(\beta_1,\ldots,\beta_n)$
with $1 \leq \beta_1< \ldots<\beta_n \leq m$,  there exists a $\C^{n+e}$-valued
holomorphic function defined on an open subset of $\mathbb{C}^{m+d}
\times J^{K}(\mathbb{C}^{m+d},\mathbb{C}^{n+e})_{(0,0)} \times
J^{K}(\mathbb{C}^{m+d},\mathbb{C}^{n+e})_{(0,0)}$ of the form
\begin{equation} \label{eq:eq2.10}
\Phi^{\alpha,\beta}(Z,\Lambda,\Gamma)=\sum_{\gamma} \frac{
R_{\gamma}^{\alpha,\beta}(\Lambda, \Gamma)}
{\big(\det(\Lambda^{r,\alpha_j})_{1\leq r,j \leq
n}\big)^{s_{\alpha\beta\gamma}} \big(
\det(\Gamma^{r,\beta_j})_{1\leq r,j \leq
n}\big)^{t_{\alpha\beta\gamma}}     } Z^{\gamma} ,
\end{equation}
where  $R^{\alpha,\beta}_\gamma$ are $\mathbb{C}^{n+e}$-valued
polynomials  and
$s_{\alpha\beta\gamma}$ and $t_{\alpha\beta\gamma}$ are nonnegative
integers, such that if $\mathcal{H}(Z,\zeta) =
\big(H(Z),\widetilde{H}(\zeta)\big)$ is an HSPM satisfying condition
$D_{\mu\nu}$, then
\begin{equation} \label{eq:eq2.11}
H(Z) = \Phi^{\mu,\nu} \Big(Z, j_0^{K}(H), j_0^K(\widetilde{H})   \Big) ,
\end{equation}
\begin{equation} \label{eq:eq2.12}
\widetilde{H}(\zeta) = \overline{\Phi^{\nu,\mu}} \Big(\zeta,
j_0^{K}(\widetilde{H}), j_0^K(H)   \Big) .
\end{equation}
Furthermore, for any $(\Lambda_0, \Gamma_0)$ such that
$\det(\Lambda^{r,\alpha_j}_0)_{1\leq r,j \leq n} \neq 0$ and
$\det(\Gamma^{r,\beta_j}_0)_{1\leq r,j \leq n}$
$ \neq 0$,
$\Phi^{\alpha,\beta}$ is holomorphic in a neighborhood of
$(0,\Lambda_0,\Gamma_0)$.
\end{finitejet}

\begin{remark3} \label{thm:NiceExample}
It is implicit in the hypotheses of Theorem \ref{thm:Finitejet} that $m \geq n$. However, if we assume that $m<n$, even if the matrices
$\big(f_z(0)\big) := \big(f^j_{z_l}(0)\big)_{1\leq l \leq m, 1  \leq j \leq n}$ and $\big(\tilde{f}_\chi(0)\big)
:= \big(\f^j_{\X_l}(0)\big)_{1\leq l \leq m, 1  \leq j \leq n}$ have maximal rank, the theorem will not hold. Let $M \subseteq \mathbb{C}^4$
be defined by $M=\{\text{Im } w_1 =|z_1|^2, \text{Im } w_2=|z_2|^2\}$. Let $M' \subseteq \mathbb{C}^4$ be defined by
$M'=\{\text{Im } w' =|z'_1|^2+|z'_2|^2+|z'_3|^2\}$.  Then $M$ is of finite type at 0, and $M'$ is 1-nondegenerate at 0.
For any positive integer $r$, define
$$ \mathcal{H}_r (z,w,\X,\T) = $$
\begin{equation}
(z_1,z_2,w_1,w_1+w_2,\chi_1-2i\chi_1 \tau_1-2i\chi_1\tau_1^r,\chi_2,\tau_1^r+\tau_1,\tau_1+\tau_2-2i\tau_1^2-2i\tau_1^{r+1}).
\end{equation}
 Observe that $\mathcal{H}_r$ is an HSPM sending $(\M,0)$ into $(\M',0)$ which is a biholomorphism  near 0.
 \end{remark3}

The proof of Theorem \ref{thm:Finitejet} will be based on arguments from
\cite{BER99b} and \cite{BRZ01}. Before proving the theorem, we first introduce
a few lemmas.

\begin{lemma2.2} \label{thm:lemma2.3}
Let $\mathcal{H}(Z,\zeta)=\big(H(Z),\tilde{H}(\zeta)\big)$ be an HSPM
sending $(\M,0)$ into $(\M',0)$.  Then
$\mathcal{H}'(Z,\zeta)=$  $\left( \bar{\tilde{H}}(Z),\bar{H}(\zeta)\right)$ is
an HSPM sending $(\M,0)$ into $(\M',0)$.
\end{lemma2.2}

\begin{proof}
Let $\rho_1, \ldots, \rho_{d}$ be defining functions for $M$, and let $\rho'_1, \ldots, \rho'_{e}$ be defining functions for $M'$.
For $j=1, \ldots, e$ and $k=1,\ldots,d$, there exist
holomorphic functions $a_k^j$ such that:
\begin{eqnarray}
\rho_j' \left( H(Z),\tilde{H}(\zeta) \right) & = & \sum_{k=1}^d a_k^j(Z,\zeta)\rho_k(Z,\zeta)  \Rightarrow    \label{q1q1}   \\
\bar{\rho}_j' \left(\tilde{H}(\zeta),H(Z)\right) & = & \sum_{k=1}^d a_k^j(Z,\zeta)\rho_k(Z,\zeta)   \Rightarrow     \label{q2q2}  \\
\rho_j' \left( \,
\bar{\tilde{H}}(Z),\bar{H}(\zeta)
\right) &  = &  \sum_{k=1}^d \bar{a}_k^j(\zeta,Z)\bar{\rho}_k(\zeta,Z)  =
\sum_{k=1}^d \bar{a}_k^j(\zeta,Z)\rho_k(Z,\zeta).    \label{q3q3}
\end{eqnarray}
 Equations (\ref{q2q2}) and (\ref{q3q3}) follow from the reality of the $\rho_j$. The result follows.
 \end{proof}

\begin{lemma2.3} \label{thm:lemma2.4}
Let $M$ and $M'$ be as in Theorem \ref{thm:Finitejet}.   Then for any
$\beta=(\beta_1,\ldots,\beta_n)$ and
$\alpha=(\alpha_1,\ldots,\alpha_n)$, with $1 \leq \alpha_1 < \ldots
< \alpha_n\leq m$, there exists a  $\C^e$-valued holomorphic function
$\phi^{\alpha}_\beta$ defined on an open subset of
$\C^{K_\beta} \times \mathbb{C}^{m+d} \times
\mathbb{C}^{m+d}$, for some integer $K_\beta$,  of the form
\begin{equation} \label{eq:eq4}
\phi^{\alpha}_\beta(\Lambda, Z,\zeta)= \sum_{\gamma, \delta}
\frac{ P_{\gamma,\delta}^{\alpha,\beta}(\Lambda) }
{\big(\det(\Lambda^{j,\alpha_l})_{1\leq j,l \leq
n}\big)^{t_{\alpha\beta\gamma\delta }}} Z^\gamma \zeta^\delta ,
\end{equation}
where  $t_{\alpha\beta\gamma\delta}$ are nonnegative integers and
 $P_{\gamma,\delta}^{\alpha,\beta}$ are $\mathbb{C}^e$-valued
polynomials, such that if $\mathcal{H}$ is an HSPM
satisfying condition $D_{\mu\nu}$, then for $(Z,\zeta) \in \mathcal{M}$,
\begin{equation} \label{eq:eq5}
Q'_{z'^\beta}\big(f(Z),\tilde{f}(\zeta),\tilde{g}(\zeta)\big) =
\phi^{\mu}_\beta \Big( j^{|\beta|}_{Z} (H),Z,\zeta \Big)  ,
\end{equation}
\begin{equation} \label{eq:eq6}
\overline{Q}'_{\chi'^\beta}\big(\tilde{f}(\zeta),{f}(Z),{g}(Z)\big) =
\overline{\phi^{\nu}_\beta} \Big( j^{|\beta|}_{\zeta} (\widetilde{H})
,\zeta,Z  \Big).
\end{equation}
%where we write $\zeta '=(\chi ',\tau ')$.
Furthermore, for  any $\Lambda_0$ such that
$\det(\Lambda_0^{j,\alpha_l})_{1 \leq j,l \leq n} \neq 0$,  $\phi_\beta ^\alpha$ is holomorphic near $(\Lambda_0,0,0)$.
\end{lemma2.3}

\begin{proof}
 For $j=1,\ldots,m$,
\begin{equation}
L_j=\frac{\partial}{\partial z_j} + \sum_{r=1}^d
Q^r_{z_j}(z,\chi,\tau) \frac{\partial}{\partial w_r}
\end{equation}
are vector fields tangent to $\mathcal{M}$. Let $\hat{z}
:=(z_{\mu_1}, \ldots, z_{\mu_n})$. Now we apply $L_{\mu_1}, \ldots,
L_{\mu_n}$ to
\begin{equation}
g(z,w) = Q' \big(f(z,w),\widetilde{H}(\chi,\tau)\big)
\end{equation}
 to get (in matrix notation):

$$g_{\hat{z}}(z,w)+   Q_{\hat{z}}(z,\chi,\tau)g_w(z,w) = $$
\begin{equation}
\Big( f_{\hat{z}}(z,w) + Q_{\hat{z}}(z,\chi,\tau) f_w(z,w) \Big)
Q'_{z'} \big(f(z,w),\widetilde{H}(\chi,\tau)\big)
\end{equation}
for all $(z,w,\chi,\tau) \in \mathcal{M}$. By assumption, $\big(
f_{\hat{z}}(0) \big)$ is invertible, so near
$(z,w,\chi,\tau)=(0,0,0,0)$, we have

$$Q'_{z'} \big(f(z,w),\widetilde{H}(\chi,\tau)\big) =$$
\begin{equation} \label{eq:eq10}
\Big( f_{\hat{z}}(z,w) + Q_{\hat{z}}(z,\chi,\tau) f_w(z,w)
\Big)^{-1} \Big( g_{\hat{z}}(z,w)+
Q_{\hat{z}}(z,\chi,\tau)g_w(z,w)\Big)  .
\end{equation}

We claim that the right hand side of (\ref{eq:eq10}) can be written in the form
\begin{equation} \label{eq:eq11}
 \sum_{\gamma,\delta}   \frac{ p_{\gamma,\delta}^{\mu} \big(j^1_{Z}(H)\big) }
   {\big(\det(f_{\hat{z}}(Z))\big)^{s_{\mu\gamma\delta }}} Z^\gamma\zeta^\delta  ,
\end{equation}
where each $p_{\gamma,\delta}^{\mu}$ is an $n \times e$ polynomial matrix and
each $s_{\mu\gamma\delta}$ is a nonnegative integer. This comes from
writing the right hand side in the following way:
\begin{equation} \label{eq:eq12}
(f_{\hat{z}}+Q_{\hat{z}}f_w)^{-1}(g_{\hat{z}}+Q_{\hat{z}}g_w) =
 (I+f_{\hat{z}}^{-1}f_wQ_{\hat{z}})^{-1} (f_{\hat{z}}^{-1})
(g_{\hat{z}}+Q_{\hat{z}}g_w)
\end{equation}
The right hand side of (\ref{eq:eq12}) has three factors. The last factor can clearly be
written in the form (\ref{eq:eq11}),
 as it is independent of $\det\big(f_{\hat{z}}(Z)\big)$. The second factor can be written in the form (\ref{eq:eq11}) since for any
 invertible matrix $A$, we can write $A^{-1}$ as $\frac{1}{\det A}(\text{adj } A)$. The first factor can also
 be written in the form (\ref{eq:eq11}). Indeed, as $f_{\hat{z}}^{-1}(0)f_w(0)Q_{\hat{z}}(0) = 0$,
  then for $(z,\chi,\tau)$ sufficiently close to 0,  $(I+B)^{-1} = \sum_{j=0}^\infty (-1)^j B^j$,
  where we define  $B:=f_{\hat{z}}^{-1}f_wQ_{\hat{z}}$.  We then use the aforementioned
 formula for the inverse of a matrix, and
 the claim is proved.

 We get (\ref{eq:eq5}) from (\ref{eq:eq10}) and (\ref{eq:eq11}) by  inductively applying the
$L_j$ and utilizing the chain rule. To complete the proof of the lemma,
we use Lemma \ref{thm:lemma2.3} to see that $\left(  \bar{\tilde{H}},\bar{H}\right)$ sends
$\mathcal{M}$ into $\mathcal{M}'$ and satisfies condition
$D_{\nu\mu}$. So as we have seen in this proof,
\begin{equation}
Q'_{z'^\beta}\left( \, \bar{\tilde{f}}(Z),\bar{f}(\zeta),\bar{g}(\zeta)\right)
= \phi^{\nu}_\beta \Big( j^{|\beta|}_{Z} \big(\bar{\tilde{H}} \big),Z,\zeta
\Big).
\end{equation}
Taking the complex conjugate of this entire equation gives  (\ref{eq:eq6}), and the proof of the lemma is complete.
\end{proof}

The following notation will be used in Lemmas \ref{thm:lemma2.5} and \ref{lemmalemma}. Let
$M$, $M'$, and $\mathcal{H}$ be as in Theorem \ref{thm:Finitejet}. We will write $j_Z^K H
=\Big( (j_Z^K)'' H, \big(g_{z^\alpha}(Z) \big)_{|\alpha| \leq K}
\Big)$, where $(j_Z^K)'' H$ represents the remaining derivatives of
$H$ at $Z$. Given any $\Lambda \in J_Z^K(\C^{m+d},\C^{n+e})_{(Z,H(Z))}$, we
will also write
\begin{equation} \label{decomposition}
\Lambda=(\Lambda_1,\Lambda_2),
\end{equation}
where $(j_Z^K H)_2$ is exactly $\big(g_{z^\alpha}(Z) \big)_{|\alpha| \leq K}$. We do a
similar decomposition for $j_\zeta^K \widetilde{H}$.

\begin{lemma2.4} \label{thm:lemma2.5}
Let $M$ and $M'$ be
as in Theorem \ref{thm:Finitejet}. For each $\beta=(\beta_1,\ldots,\beta_{n+e})$  and
$\alpha=(\alpha_1,\ldots,\alpha_n)$ with $1\leq \alpha_1 < \ldots <
\alpha_n \leq m$,
 there exists  a $\C^{n+e}$-valued function $\Psi_\beta^\alpha$, holomorphic on an open subset of
 $\mathbb{C}^{m+d} \times \mathbb{C}^{m+d} \times  \C^{K_\beta}$
 for some integer $K_\beta$, of the form
\begin{equation} \label{eq:eq14}
\Psi_\beta ^\alpha (Z,\zeta,\Lambda)=\sum_{\gamma, \delta, \kappa}
\frac {P^{\alpha,\beta}_{\gamma,\delta,\kappa}(\Lambda_1)}
{\big(\det(\Lambda^{r,\alpha_l})_{1 \leq l,r \leq n}\big)^{t_{\alpha
\beta \gamma \delta\kappa}}} Z^\gamma \zeta^\delta\Lambda_2^\kappa ,
\end{equation}
where   $P_{\gamma,\delta,\kappa}^{\alpha,\beta}(\Lambda_1)$ are
$\mathbb{C}^{n+e}$-valued polynomials  and  $t_{\alpha
\beta\gamma\delta\kappa}$ are nonnegative integers,     such that
if $\mathcal{H}$ is an HSPM satisfying condition $D_{\mu\nu}$, then for
$(Z,\zeta) \in \mathcal{M}$,
\begin{equation} \label{eq:eq15}
\partial^\beta H(Z) = \Psi_\beta ^\nu \big(Z,\zeta, j^{k+|\beta|}_{\zeta}(\widetilde{H})        \big),
\end{equation}

\begin{equation} \label{eq:eq16}
\partial^\beta \widetilde{H}(\zeta) = \overline{\Psi_\beta ^\mu} \big(\zeta,Z, j^{k+|\beta|}_{Z}(H)       \big).
\end{equation}
Furthermore, for  any $\Lambda_0$ such that
$\det(\Lambda_0^{r,\alpha_l})_{1 \leq l,r \leq n} \neq 0$, $\Psi_\beta ^\alpha$ is holomorphic near $(0,0,\Lambda_0)$.
\end{lemma2.4}

\begin{proof}
 As $M'$ is $k$-nondegenerate at 0, assume the vectors
$\overline{Q}'^{j_1}_{z\chi^{\alpha_1}}(0),\ldots,
\overline{Q}'^{j_n}_{z\chi^{\alpha_n}}(0)$ span $\mathbb{C}^n$ where each
$j_k \in \{1,\ldots,e\}$, each $|\alpha_j| \leq k$, and
$\overline{Q}'=(\overline{Q}'^1,\ldots, \overline{Q}'^e)$.
%and $\overline{Q}'^j_{z\X^\alpha} = (\overline{Q}'^j_{z_1\X^\alpha},\ldots,$ $\overline{Q}'^j_{z_m\X^\alpha})$. From Lemma
\ref{thm:lemma2.4}, we have for each $(Z,\zeta) \in \M$:
\begin{eqnarray}    \label{eq2.40z}
\overline{Q}'^{j_1}_{\chi'^{\alpha_1}}(\tilde{f}(\zeta),f(Z),g(Z)) & = &\overline{(\phi^{\nu}_{\alpha_1})^{j_1}}
\Big( j^{|\alpha_1|}_{\zeta} (\widetilde{H}) ,\zeta,Z  \Big) \nonumber \\
\vdots \,\,\,\,\,\,\, & \vdots & \,\,\,\,\,\,\, \vdots \nonumber \\
\overline{Q}'^{j_n}_{\chi'^{\alpha_n}}(\tilde{f}(\zeta),f(Z),g(Z)) & = &
\overline{(\phi^{\nu}_{\alpha_n})^{j_n}} \Big(
j^{|\alpha_n|}_{\zeta} (\widetilde{H}) ,\zeta,Z  \Big) ,
\end{eqnarray}
where $\phi_\alpha^\beta = \left( (\phi_\alpha^\beta)^1, \ldots, (\phi_\alpha^\beta)^e \right)$.
Using this system of equations, coupled with the fact that normal
coordinates for $M'$ imply that $\overline{Q}'(\chi',0,w')
\equiv\overline{Q}'(0,z',w') \equiv w'$, we can apply the implicit function
theorem to find a map $B^\nu: \mathbb{C}^n \times \mathbb{C}^n
\rightarrow \mathbb{C}^{n+e}$, holomorphic
near 0, such that
\begin{equation} \label{eq:eq18}
H(Z) = B^\nu \Bigg( \tilde{f}(\zeta) , \bigg(
\overline{(\phi^{\nu}_{\alpha_l})^{j_l}} \Big( j^{|\alpha_l|}_{\zeta}
\widetilde{H} ,\zeta,Z\Big)  \bigg)_{1 \leq l \leq n} \Bigg).
\end{equation}
Now we are going to write each $
\overline{(\phi^{\nu}_{\alpha_l})^{j_l}} \left( j^{|\alpha_l|}_{\zeta}
\widetilde{H} ,\zeta,Z\right)$ in a different form. For $1 \leq j \leq m$,
the following vector fields are tangent to $\mathcal{M}$:
\begin{equation}
 \displaystyle{\widetilde{L}_j = \frac{\partial}{\partial \chi_j}+
\sum_{r=1}^d \overline{Q}^r_{\chi_j}(\chi,z,w)\frac{\partial}{\partial \tau_r}}.
\end{equation}
Apply $\widetilde{L}_{\nu_1},\ldots, \widetilde{L}_{\nu_n}$ to
$\tilde{g}(\chi,\tau)=\overline{Q}'\big(\tilde{f}(\chi,\tau),f(z,w),g(z,w)\big)$
repeatedly ($|\beta|$ times), and apply Cramer's rule each time to
see that for $(z,w,\chi,\tau) \in \mathcal{M}$ and each
$l=1,\ldots,e$,
\begin{equation}
\overline{Q}'^l_{\chi'^\beta}\big(\tilde{f}(\chi,\tau),f(z,w),g(z,w)\big)=  \sum_{1
\leq |\gamma| \leq |\beta|}
\frac{\big(\widetilde{L}^\gamma\tilde{g}^l(\chi,\tau)\big)P^{\beta,l}_\gamma \Big(
\big(\widetilde{L}^\delta\tilde{f}(\chi,\tau)\big)_{1 \leq |\delta| \leq |\beta|} \Big)}
{\det \big(\widetilde{L}_{\nu_i}\tilde{f}^j(\chi,\tau) \big)_{1 \leq i,j \leq n}}   ,
\end{equation}
where each $P^{\beta,l}_\gamma$ is a polynomial independent of
$M,M', $ and $\mathcal{H}$. Notice that by assumption, the
denominator is nonzero near $(\chi,\tau)=(0,0)$. So we have
\begin{equation}
\overline{\big(\phi^{\nu}_\beta \big)^l} \big( j^{|\beta|}_{\zeta}
\widetilde{H} ,\zeta,Z  \big)= \sum_{1
\leq |\gamma| \leq |\beta|}
\frac{\big(\widetilde{L}^\gamma\tilde{g}^l(\chi,\tau)\big)P^{\beta,l}_\gamma \Big(
\big(\widetilde{L}^\delta\tilde{f}(\chi,\tau)\big)_{1 \leq |\delta| \leq |\beta|} \Big)}
{\det \big(\widetilde{L}_{\nu_i}\tilde{f}^j(\chi,\tau) \big)_{1 \leq i,j \leq n}}.
\end{equation}
Substituting this in (\ref{eq:eq18}), we get
%$$H(Z) =$$
\begin{equation} \label{eq:eq22}
H(Z)=B^\nu \left( \tilde{f}(\zeta) , \left( \sum_{1
\leq |\gamma| \leq |\alpha_l|}
\frac{\big(\widetilde{L}^\gamma\tilde{g}^{j_l}(\zeta)\big)P^{\alpha_l,j_l}_\gamma \Big(
\big(\widetilde{L}^\delta\tilde{f}(\zeta)\big)_{1 \leq |\delta| \leq |\alpha_l|} \Big)}
{\det \big(\widetilde{L}_{\nu_i}\tilde{f}^j(\zeta) \big)_{1 \leq i,j \leq n}}
\right)_{1 \leq l \leq n} \right).
\end{equation}
If we Taylor expand, we can write the right hand side of (\ref{eq:eq22}) as
\begin{equation}
\sum_{\alpha,\beta,\gamma}A^\nu_{\alpha\beta\gamma}(\widetilde{\Lambda}_1)Z^\alpha
\zeta^\beta \widetilde{\Lambda}_2^\gamma,
\end{equation}
where we remind the reader that $\widetilde{\Lambda}_2$ corresponds to
$\big(\tilde{g}_{\chi^\alpha}(\zeta)\big)$, and $\widetilde{\Lambda}_1$
corresponds to the remaining derivatives of $\widetilde{H}$ at $\zeta$. We claim that each
$A^\nu_{\alpha\beta\gamma}$ is rational. This follows from the fact
that
\begin{equation}
\widetilde{L}^\gamma\tilde{g}(\chi,\tau) =
\tilde{g}_{\chi^\gamma}(\chi,\tau) + R
\Big(j^{|\gamma|}_{\zeta} (\tilde{g}),\chi,z,w\Big),
\end{equation}
where $R$ is a holomorphic mapping which vanishes when $\chi=z=w=0$.

 Furthermore, each
$A^\nu_{\alpha\beta\gamma}(\widetilde{\Lambda}_1)$ is of the form given in the right hand side of
 (\ref{eq:eq14}). This can be seen by Taylor expanding $B^\nu$ as given in
(\ref{eq:eq18}) and plugging in (\ref{eq:eq4}). Define
\begin{equation}
\Psi_0^\nu(Z,\zeta,\widetilde{\Lambda}) :=
\sum_{\alpha,\beta,\gamma}A^\nu_{\alpha\beta\gamma}(\widetilde{\Lambda}_1)Z^\alpha
\zeta^\beta \widetilde{\Lambda}_2^\gamma  .
\end{equation}
This proves (\ref{eq:eq15}) for $|\beta|=0$. For $|\beta|>0$, as every
point in $\mathcal{M}$ is of the form
$\big(z,w,\chi,\overline{Q}(\chi,z,w)\big)$, we have that
\begin{equation} \label{qweewq}
H(z,w) \equiv \Psi^\nu_0
\left(z,w,\chi,\overline{Q}(\chi,z,w),j^k_{(\chi,\overline{Q}(\chi,z,w))}\big(\widetilde{H}\big)
\right) .
\end{equation}
We inductively differentiate (\ref{qweewq}), applying the chain rule, and (\ref{eq:eq15}) follows.

\indent To get (\ref{eq:eq16}), we know from Lemma \ref{thm:lemma2.3} that
$\left(\bar{\tilde{H}},\bar{H}\right)$ sends $\mathcal{M}$ into $\mathcal{M}'$
and satisfies condition $D_{\nu\mu}$. So
\begin{equation}
\partial^\beta \bar{\tilde{H}}(Z) = \Psi_\beta ^\mu \left(Z,\zeta, j^{k+|\beta|}_{\zeta}\big(  \bar{H}\big)        \right)  .
\end{equation}
Take the complex conjugate of both sides of this equation, and the lemma follows.
\end{proof}

Now we define the \emph{r-th Segre mappings} of $M$ at 0. These maps
were first introduced by Baouendi, Ebenfelt, and Rothschild in
\cite{BER96} and will prove extremely useful in completing the proof of
Theorem \ref{thm:Finitejet}. Given a positive integer $r$, let $t^0,\ldots,t^{r-1} \in \C^m$ and define $v^r:
\mathbb{C}^{rm} \rightarrow \mathbb{C}^{m+d}$ in the following way:
\begin{equation} \label{segreMapping1}
v^r(t^0, \ldots, t^{r-1}) := \left(t^0, u^r(t^0, \ldots,
t^{r-1})\right),
\end{equation}
where $u^r: \mathbb{C}^{rm} \rightarrow \mathbb{C}^d$ is given
inductively by
\begin{equation} \label{segreMapping2}
u^1(t^0)=0 \,\, , \,\,
u^r(t^0,\ldots,t^{r-1})=Q\left(t^0,t^1,\overline{u^{r-1}}(t^1,\ldots,
t^{r-1})\right) \text{ for } r\geq 2 .
\end{equation}

\begin{defn2.5}
Let $V$ and $W$ be finite dimensional complex vector spaces. Let
$\mathcal{R}_0(V \times W, V)$ denote the ring of germs of
holomorphic functions $f$ at $V \times \{0\}$ in $V \times W$ which
can be written in the form $f(\Lambda,\Gamma) =\sum_{\alpha}
p_\alpha(\Lambda) \Gamma^\alpha$, where each $p_\alpha(\Lambda)$ is
a polynomial function on $V$.
\end{defn2.5}

The following lemma is proved in \cite{BRZ01}:
\begin{lemma2.6} \label{thm:lemma2.7}
Let $V_0,V_1,\widetilde{V}_0,\widetilde{V}_1$ be finite dimensional complex
vector spaces with fixed bases and $x_0,x_1,\tilde{x}_0,\tilde{x}_1$
be the linear coordinates with respect to these bases. Let $p \in
\mathbb{C}[x_0]$ and $\tilde{p} \in \mathbb{C}[\tilde{x}_0]$ be
nontrivial polynomial functions on $V_0$ and $\widetilde{V}_0$
respectively, and let
$$\phi=(\phi_0,\phi_1):\mathbb{C} \times V_0 \times V_1 \rightarrow \widetilde{V}_0 \times \widetilde{V}_1$$
 be a germ of a holomorphic map with components in $\mathcal{R}_0(\mathbb{C} \times V_0 \times V_1, \mathbb{C} \times V_0)$,
  such that $\phi(\mathbb{C} \times V_0 \times \{0\}) \subseteq \widetilde{V}_0 \times \{0\}$,
  and satisfying $\displaystyle{  \tilde{p}\left( \phi_0 \left( \frac{1}{p(x_0)},x_0,0 \right)\right) \not\equiv 0}$.
  Then given any $\tilde{h} \in \mathcal{R}_0(\mathbb{C} \times \widetilde{V}_0 \times \widetilde{V}_1, \mathbb{C} \times \widetilde{V}_0)$,
  there exists $h \in \mathcal{R}_0(\mathbb{C} \times V_0 \times V_1, \mathbb{C} \times V_0)$ such that
\begin{equation}
\tilde{h}\left( \frac{1}{\tilde{p}\left(
\phi_0\left(\frac{1}{p(x_0)},x_0,x_1\right)\right)},\phi \left(
\frac{1}{p(x_0)},x_0,x_1\right)\right) \equiv h \left(
\frac{1}{q(x_0)},x_0,x_1\right),
\end{equation}
with $q(x_0):=p(x_0)^t \tilde{p}\left(
\phi_0\left(\frac{1}{p(x_0)},x_0,0\right)\right)$ for some positive integer $t$.
Furthermore, $h$ vanishes on $\mathbb{C} \times V_0 \times \{0\}$ if
$\tilde{h}$ vanishes on $\mathbb{C} \times \widetilde{V}_0 \times
\{0\}$.
\end{lemma2.6}

Lemma \ref{thm:lemma2.7} will be key in establishing the following lemma.
\begin{lemma2.7}    \label{lemmalemma}
Let $M$ and $M'$ be as in Theorem \ref{thm:Finitejet}.
Given any $\beta=(\beta_1,\ldots,\beta_n)$ and
$\alpha=(\alpha_1,\ldots,\alpha_n)$, with $1 \leq \beta_1 < \ldots <
\beta_n\leq m$ and $1 \leq \alpha_1 < \ldots < \alpha_n\leq m$, and
any positive integer $s$, there exists a $\C^{n+e}$-valued function
$\Xi^{\alpha,\beta}_s(x,\Lambda,\Gamma)$ holomorphic on an open
subset of $\mathbb{C}^{sm} \times J^{sk}
(\mathbb{C}^{m+d},\mathbb{C}^{n+e})_{(0,0)} \times J^{sk}
(\mathbb{C}^{m+d},\mathbb{C}^{n+e})_{(0,0)}$ of the form
\begin{equation} \label{eq:niceForm}
\Xi^{\alpha,\beta}_s(x,\Lambda,\Gamma) = \sum_{\gamma}
\frac{P^{\alpha\beta s}_{\gamma}(\Lambda,\Gamma)}{Q^{\alpha\beta
s}_{\gamma}(\Lambda,\Gamma)}x^\gamma ,
\end{equation}
where each $P^{\alpha\beta s}_{\gamma}$ is a $\C^{n+e}$-valued
polynomial, and
$$Q^{\alpha\beta
s}_{\gamma}(\Lambda,\Gamma):=\big(\det(\Lambda^{r,\alpha_l})_{1 \leq
r,l \leq n}\big)^{u_{\alpha\beta\gamma s}}
\big(\det(\Gamma^{r,\beta_l})_{1 \leq r,l \leq
n}\big)^{v_{\alpha\beta\gamma s}}$$ for some
 nonnegative integers
$u_{\alpha\beta\gamma s}$ and $v_{\alpha\beta\gamma s}$, such that
if $\mathcal{H}$ is an HSPM satisfying condition $D_{\mu\nu}$, then
\begin{equation} \label{eq:eq32}
H\big( v^s(t^0,\ldots,t^{s-1})\big) =
\Xi^{\mu,\nu}_s\big(t^0,\ldots,t^{s-1},j_0^{sk}H, j_0^{sk}\widetilde{H}
\big) ,
\end{equation}
\begin{equation} \label{eq:eq33}
\,\,\,\,\, \widetilde{H}\big( \overline{v^s}(t^0,\ldots,t^{s-1})\big) =\overline{
\Xi^{\nu,\mu}_s} \big(t^0,\ldots,t^{s-1},j_0^{sk}\widetilde{H}, j_0^{sk}
{H} \big) \,\,\, .
\end{equation}
Furthermore, for any $(\Lambda_0, \Gamma_0)$ such that
$\det(\Lambda^{r,\alpha_l}_0)_{1\leq r,l \leq n} \neq 0$ and
$\det(\Gamma^{r,\beta_l}_0)_{1\leq r,l \leq n} \neq 0$,
$\Xi^{\alpha,\beta}_s$ is holomorphic on a neighborhood of
$(0,\Lambda_0,\Gamma_0)$.
\end{lemma2.7}

\begin{proof}
We inductively prove
something stronger. First, we simplify notation slightly. Define
$$p^\alpha\big(\Lambda_1(Z)\big) :=\det\left(
\Lambda^{r,\alpha_l}(Z)\right)_{1 \leq r,l \leq n} \, ,$$
\begin{equation}
\tilde{p}^\beta \big(\widetilde{\Lambda}_1 (\zeta)\big) :=
\det( \widetilde{\Lambda}^{r,\beta_l}(\zeta))_{1 \leq r,l \leq n} \, ,
\end{equation}
where $\Lambda_1$  is as defined in (\ref{decomposition}) ($\widetilde{\Lambda}_1$ is defined in a similar way).
 We will show that for any
$\gamma$ and $s$, there exist nonnegative integers $a_{\alpha\beta}^s$ and $b_{\alpha\beta}^s$
 and  holomorphic maps $\Theta^{\alpha,\beta,\gamma}_s$ with components in
$$\mathcal{R}_0\Big(\mathbb{C} \times J^{ks+|\gamma|}(\C^{m+d},\C^{n+e})_{(0,0)} \times
J^{ks+|\gamma|}
 (\C^{m+d},\C^{n+e})_{(0,0)} \times \mathbb{C}^{ms},$$
$$ \mathbb{C} \times
J^{ks+|\gamma|}(\C^{m+d},\C^{n+e})_{(0,0)}  \times J^{ks+|\gamma|} (\C^{m+d}, \C^{n+e})_{(0,0)} \Big)$$ such that
$$\partial^\gamma H\big( v^s(t^0,\ldots,t^{s-1})\big) =$$
\begin{equation} \label{eq:eq34}
  \Theta^{\mu,\nu,\gamma}_s
\left(\frac{1}{p^\mu \big(\Lambda_1(0)\big)^{a_{\mu\nu}^s} \tilde{p}^\nu
\big( \widetilde{\Lambda}_1(0)\big)^{b_{\mu\nu}^s}} , j^{ks+|\gamma|}_0 H, j^{ks+|\gamma|}_0 \widetilde{H},t^0,\ldots,t^{s-1} \right),
\end{equation}

$$\partial^\gamma \widetilde{H}\big( \overline{v^s}(t^0,\ldots,t^{s-1})\big) =$$
\begin{equation} \label{eq:eq35}
 \overline{\Theta^{\nu,\mu,\gamma}_s}
\left(\frac{1}{{p}^\nu \big(\widetilde{\Lambda}_1(0)\big)^{a_{\nu\mu}^s} \tilde{p}^\mu \big( {\Lambda}_1(0)\big)^{b_{\nu\mu}^s}} ,
j^{ks+|\gamma|}_0 \widetilde{H}, j^{ks+|\gamma|}_0 {H},t^0,\ldots,t^{s-1} \right).
\end{equation}

 First, we will use new
notation to reformulate Lemma \ref{thm:lemma2.5}. Let $j_Z^l H= (\hat{j}_Z^l
H, \hat{\hat{j}}_Z^l H)$, where $\hat{\hat{j}}_Z^l H = \big(
g_{z^\alpha}(Z)\big)_{|\alpha| \leq l}$, and $\hat{j}_Z^l H$
represents the remaining derivatives at $Z$. (A similar decomposition applies
to  $\widetilde{H}$).  According to Lemma \ref{thm:lemma2.5}, there exist maps
$\theta^\alpha_\gamma $ with components in $\mathcal{R}_0(\mathbb{C} \times
\mathbb{C}^{l'_{\gamma}} \times \mathbb{C}^{2m+2d} \times
\mathbb{C}^{l''_{\gamma}},  \mathbb{C} \times
\mathbb{C}^{l'_{\gamma}} )$ (for some integers $l'_{\gamma}$ and
$l''_{\gamma}$) such that for $(Z,\zeta) \in \mathcal{M}$:
\begin{equation} \label{eq:eq36}
\partial^\gamma H(Z) = \theta^\nu_\gamma \Bigg( \frac{1}{\tilde{p}^\nu(\hat{j}_\zeta^{k+|\gamma|} \widetilde{H})},
\hat{j}_\zeta^{k+|\gamma|} \widetilde{H},Z,\zeta, \hat{\hat{j}}_\zeta^{k+|\gamma|} \widetilde{H} \Bigg) ,
\end{equation}
\begin{equation} \label{eq:eq37}
\partial^\gamma \widetilde{H}(\zeta) = \overline{\theta^\mu_\gamma} \Bigg( \frac{1}{{p}^\mu(\hat{j}_Z^{k+|\gamma|} {H})},
\hat{j}_Z^{k+|\gamma|}{H},
 \zeta, Z,  \hat{\hat{j}}_Z^{k+|\gamma|} {H} \Bigg).
\end{equation}
It is easy to show that (\ref{eq:eq34}) and (\ref{eq:eq35}) hold for $s=1$ by letting $(Z,\zeta)=\big( (z,0),0\big)$ in (\ref{eq:eq36}) and
$(Z,\zeta)=\big(0,(\chi,0)\big)$ in (\ref{eq:eq37}). So now assume for some $s>1$, (\ref{eq:eq34}) and (\ref{eq:eq35}) hold for $s-1$. We will
show they hold for $s$.

For any $s$, it is clear from the definition of the Segre mappings that
\begin{equation}
\big(
v^s(t^0,\ldots,t^{s-1}), \overline{v^{s-1}}(t^1,\ldots,t^{s-1})
\big) \in \mathcal{M}.
\end{equation}
 Using this fact in (\ref{eq:eq36}), we see that
$$\partial^\gamma H\big( v^s(t^0,\ldots,t^{s-1})\big)  =
\theta^\nu_\gamma \left(
\frac{1}{\tilde{p}^\nu\Big(\hat{j}_{\overline{v^{s-1}}(t^1,\ldots,t^{s-1})}^{k+|\gamma|}
\widetilde{H}\Big)}, \right.$$
\begin{equation} \label{eq:eq38}
   \hat{j}_{\overline{v^{s-1}}(t^1,\ldots,t^{s-1})}^{k+|\gamma|} \widetilde{H}, v^s(t^0,\ldots,t^{s-1}), \overline{v^{s-1}}(t^1,\ldots,t^{s-1}),
  \hat{\hat{j}}_{\overline{v^{s-1}}(t^1,\ldots,t^{s-1})}^{k+|\gamma|} \widetilde{H} \Bigg).
\end{equation}
But by our induction hypothesis,
\begin{equation*}
 j_{\overline{v^{s-1}}(t^1,\ldots,t^{s-1})}^{k+|\gamma|} \widetilde{H} =
\end{equation*}
\begin{equation}   \label{eq:eq39}
\Bigg(      \overline{\Theta^{\nu,\mu,\Delta}_{s-1}}
\Bigg(\frac{1}{\tilde{p}^\nu
\big(\widetilde{\Lambda}_1(0)\big)^{a_{\nu\mu}^{s-1}} p^\mu \big(
{\Lambda}_1(0)\big)^{b_{\nu\mu}^{s-1}}} ,
j_0^{ks+|\Delta|-k} {H}, j^{ks+|\Delta|-k}_0 \widetilde{H}, t^1,\ldots,t^{s-1} \Bigg)\Bigg)_{|\Delta|
\leq k+|\gamma|} .
\end{equation}
For convenience, we write the tuple on the right hand side of (\ref{eq:eq39}) as $(A,B)$
where $B$ corresponds to
 $\Big( \tilde{g}_{\chi^\alpha}\big(\overline{v^{s-1}}(t^1,\ldots,t^{s-1})\big)\Big)_{|\alpha| \leq k+|\gamma|}$,
 and $A$ corresponds to the remainder.
We plug (\ref{eq:eq39}) into (\ref{eq:eq38}) to get
\begin{equation}
\partial^\gamma H\big( v^s(t^0,\ldots,t^{s-1})\big) = \theta^\nu_\gamma
\Bigg( \frac{1}{\tilde{p}^\nu(A)},A,v^s(t^0,\ldots,t^{s-1}),\overline{v^{s-1}}(t^1,\ldots,t^{s-1}),B \Bigg) .
\end{equation}
Thus, (\ref{eq:eq34}) follows from Lemma \ref{thm:lemma2.7}.
%The correspondences between the lemma and the notation above are as follows:
%\begin{eqnarray}
%p & \sim & \tilde{p}^\nu(\widetilde{\Lambda}(0))^{a_{\nu\mu}^{s-1}} p^\mu(\Lambda(0))^{b_{\nu\mu}^{s-1}} \nonumber \\
%\tilde{p} & \sim & \tilde{p}^\nu(\widetilde{\Lambda}(0)) \nonumber \\
% x_0 & \sim & \big(j_0^{ks+|\gamma|} \widetilde{H}, j_0^{ks+|\gamma|} H \big) \nonumber \\
% x_1 & \sim & (t^0,\ldots,t^{s-1}) \nonumber \\
%  \phi & \sim & \big( A,v^s(t^0,\ldots,t^{s-1}),\overline{v^{s-1}}(t^1,\ldots,t^{s-1}),B \big) \nonumber \\
%\phi_0 & \sim & A \nonumber \\
%  \tilde{h} & \sim  & \theta_\gamma^\nu
%  \end{eqnarray}

 To finish the proof, we need only show (\ref{eq:eq35}). Here we apply
Lemma \ref{thm:lemma2.3}, which tells us that $\left(\bar{\tilde{H}},\bar{H}\right)$ sends
$\mathcal{M}$ into $\mathcal{M}'$ and satisfies condition
$D_{\nu\mu}$. So by (\ref{eq:eq34}), we see that

$$\partial^\gamma \bar{\tilde{H}}\big( v^s(t^0,\ldots,t^{s-1})\big) = $$
\begin{equation} \label{eq:eq41}
\Theta^{\nu,\mu,\gamma}_s \Bigg(\frac{1}{p^\nu \big(\bar{\tilde{\Lambda}}_1(0)\big)^{a_{\nu\mu}^s} \tilde{p}^\mu
 \big( \bar{\Lambda}_1(0)\big)^{b_{\nu\mu}^s}} , j^{ks+|\gamma|}_0 \bar{\tilde{H}}, j^{ks+|\gamma|}_0 \bar{H}, t^0,\ldots,t^{s-1} \Bigg) .
\end{equation}
As $p^\nu$ and $\tilde{p}^\mu$ are polynomials with real coefficients, we take the complex conjugate of both sides of (\ref{eq:eq41}) to see
that (\ref{eq:eq35}) holds true.
\end{proof}

We are  almost ready to complete the proof of Theorem \ref{thm:Finitejet}. First,
however, we present  three lemmas. Lemma \ref{thm:lemma2.9} can be
found (using slightly different language) in \cite{BER99b} and is thus presented here without proof.  Lemma
\ref{thm:lemma2.10}
is a generalization of a lemma  found in \cite{BRZ01}.
Lemma \ref{thm:lemma2.11} can be found in \cite{BRZ01} and is  presented here without proof.

\begin{lemma2.8} \label{thm:lemma2.9}
Let $M$ be as in Theorem \ref{thm:Finitejet}. Then there exists an integer $r$ such
that the matrix
\begin{equation} \label{eq:eq42}
\left( \frac{\partial v^{2r} }{\partial
(t^0,t^{r+1},t^{r+2},\ldots,t^{2r-1})}
(0,x^1,\ldots,x^{r-1},x^r,x^{r-1},\ldots,x^1) \right)
\end{equation}
has rank $m+d$ for all $(x^1,\ldots,x^r) \in U \backslash V$,
where $U \subseteq \C^{rm}$ is an open neighborhood of the origin, and $V$ is a proper
holomorphic subvariety of $U$. In
addition,
\begin{equation} \label{eq:eq43}
v^{2r} (0,x^1,\ldots,x^{r-1},x^r,x^{r-1},\ldots,x^1) \equiv 0 .
\end{equation}
(Here, $v^{2r}$ is as defined in (\ref{segreMapping1}).)
\end{lemma2.8}

\begin{lemma2.9} \label{thm:lemma2.10}
Let $V :(\mathbb{C}^{r_1} \times \mathbb{C}^{r_2},0) \rightarrow
(\mathbb{C}^N,0)$, $r_2 \geq N$, be a holomorphic map, defined near $0$,
satisfying $V(x,\xi) \big|_{\xi=0} \equiv 0$, with $(x,\xi) \in
\mathbb{C}^{r_1} \times \mathbb{C}^{r_2}$, and assume the matrix
 $ \big( \frac{\partial V}{\partial \xi}(x,0)\big)$ has an
 $N \times N$ minor which is not identically $0$.
Then there exist holomorphic maps (defined near $0$)
\begin{equation} \label{eq:eq44}
\delta:(\mathbb{C}^{r_1},0) \rightarrow \mathbb{C} \,\, , \,\,
\phi:(\mathbb{C}^{r_1} \times \mathbb{C}^N ,0) \rightarrow
(\mathbb{C}^{r_2},0),
\end{equation}
with $\delta(x) \not\equiv 0$ such that
\begin{equation} \label{eq:eq45}
V \left( x, \phi \left( x, \frac{Z}{\delta(x)}\right) \right) \equiv Z
\end{equation}
for all $(x,Z) \in \mathbb{C}^{r_1} \times \mathbb{C}^N$  such that
$\delta(x) \neq 0$ and both $x$ and $\frac{Z}{\delta(x)}$ are
sufficiently small.  Furthermore, if $V$ is holomorphic algebraic, then
given any sufficiently small $x_0$  satisfying $\delta(x_0) \neq 0$,
the map $\varphi_{x_0}(Z) := \phi\left( x_0, \frac{Z}{\delta(x_0)}\right)$ is
holomorphic algebraic for all $Z$ in a neighborhood of $0$.
\end{lemma2.9}

\begin{proof}
Write $\xi = (\xi',\xi'')$, where $\xi'=(\xi_1,\ldots,\xi_N) \in \C^N$ and $\xi''
=(\xi_{N+1},\ldots,\xi_{r_2}) \in \C^{r_2-N}$.  Assume, without loss of generality, that
$\det\left(\frac{\partial V}{\partial \xi'}(x,0)\right) \not\equiv 0$.
We
wish to solve the equation
\begin{equation}  \label{solve}
Z = V(x,\xi',0)
\end{equation}
for $\xi'$.
As $V(x,0) \equiv 0$, we can write
\begin{equation}
Z = V(x,\xi',0) = a(x,\xi')\xi' ,
\end{equation}
where $a(x,\xi')$ is an $N \times N$ matrix of holomorphic functions defined near 0.
Furthermore, by expanding $a(x,\xi')$, we can write
\begin{equation}  \label{qaz}
Z = V(x,\xi',0) = a(x,0)\xi' + \left(  (\xi')^T R_j(x,\xi')\xi'\right)_{1 \leq j \leq N},
\end{equation}
where each $R_j(x,\xi')$ is an $N \times N$ matrix of holomorphic functions defined near 0.
Define $d(x) := \det\left(\frac{\partial V}{\partial \xi'}(x,0)\right)$. Using the fact that
$\big(\text{adj}(A)\big)A = \det(A)I$ for any square matrix $A$, we multiply the far left and far right
sides of (\ref{qaz}) by $b(x):=\text{adj}\big(a(x,0)\big)$, noting that $a(x,0) = \frac{\partial V}{\partial \xi'}(x,0)$,
to get
\begin{equation}  \label{qaz1}
b(x)Z - d(x)\xi' - b(x)\left(  (\xi')^T R_j(x,\xi')\xi'\right)_{1 \leq j \leq N} = 0 .
\end{equation}
Divide both sides of (\ref{qaz1}) by $d(x)^2$, and substitute $\tilde{\xi}' = \frac{\xi'}{d(x)}$ and
$\widetilde{Z} = \frac{Z}{d(x)^2}$ to get
\begin{equation}  \label{qaz2}
b(x)\widetilde{Z} - \tilde{\xi}' - b(x)\left(  (\tilde{\xi}')^T R_j(x,d(x)\tilde{\xi}')\tilde{\xi}'\right)_{1 \leq j \leq N} = 0.
\end{equation}
By the implicit function theorem, there is a unique holomorphic solution $\tilde{\xi}' = \theta(x,\widetilde{Z})$
defined near 0 such that $\theta(0)=0$.
Thus, the first part of the  theorem follows by letting $\delta(x) := d(x)^2$ and
\begin{equation}
\phi\left(x,y\right) := \Big( d(x)
\theta\left(x,y\right) , 0 \Big) ,
\end{equation}
for $(x,y) \in \C^{r_1} \times \C^N$.
If $V$ is algebraic,  the last
part of the theorem then follows from the \emph{algebraic} implicit function theorem (see, e.g., \cite{BER99a}).
\end{proof}

\begin{lemma2.10} \label{thm:lemma2.11}
Let $V_0$ and $V_1$ be finite dimensional vector spaces with fixed linear
coordinates $x_0$ and $x_1$, respectively. Let $P(x_0,x_1,\lambda)
\in \mathcal{R}_0(V_0 \times V_1 \times \mathbb{C}, V_0)$ with
$P(x_0,0,0) \equiv 0$. For a given integer $l \geq 0$, consider the
Laurent series expansion
\begin{equation} \label{eq:eq46}
P\Big(x_0,\frac{x_1}{\lambda^l},\lambda \Big) = \sum_{\nu \in
\mathbb{Z}} c_\nu (x_0,x_1) \lambda^\nu \,\, .
\end{equation}
Then $c_0(x_0,0) \equiv 0$, and for every $\nu \in \mathbb{Z}$,
$c_\nu \in \mathcal{R}_0(V_0 \times V_1, V_0)$.
\end{lemma2.10}

%$\,$ \vspace{-.2in}
%\begin{center}
%\textbf{Proof of Theorem \ref{thm:Finitejet}}
%\end{center}

\subsection{Proof of Theorem \ref{thm:Finitejet}}

Now we prove Theorem \ref{thm:Finitejet}.
Let $r$ be as in Lemma \ref{thm:lemma2.9}. We take
$x=(x^1,\ldots,x^r) \in \mathbb{C}^{rm}$ and $y=(y^0,\ldots,y^{r-1})$
$\in \mathbb{C}^{rm}$. Let $L(x,y) :=(y^0,x^1,\ldots,x^r,
x^{r-1}+y^{r-1}, \ldots, x^1+y^1)$ and $V(x,y):= v^{2r}(L(x,y))$. In
Lemma \ref{thm:lemma2.10}, we take $r_1=r_2=rm$.  From (\ref{eq:eq43}), we see that $V(x,0)
\equiv 0$. Also, from Lemma \ref{thm:lemma2.9} we see that the other hypothesis of
Lemma \ref{thm:lemma2.10} holds. Thus we apply Lemma \ref{thm:lemma2.10}. Let $\delta$ and $\phi$ be as
given in the lemma. We plug these into (\ref{eq:eq32}) to see that
\begin{equation} \label{eq:eq56}
H(Z) \equiv \Xi_{2r}^{\mu,\nu} \left( L\left(x,\phi
\left(x,\frac{Z}{\delta(x)}\right)\right), j_0^{2rk} H, j_0^{2rk}
\widetilde{H} \right) .
\end{equation}
We rewrite the right hand side of (\ref{eq:eq56}) in the following way:
\begin{equation} \label{eq:eq57}
H(Z) \equiv \hat{\Xi}_{2r}^{\mu,\nu} \left(   j_0^{2rk} H, j_0^{2rk}
\widetilde{H},  \frac{Z}{\delta(x)}, x  \right),
\end{equation}
noting that the components of $$\hat{\Xi}_{2r}^{\mu,\nu}: J^{2rk}(\mathbb{C}^{m+d},\mathbb{C}^{n+e})_{(0,0)}
\times J^{2rk}(\mathbb{C}^{m+d},\mathbb{C}^{n+e})_{(0,0)} \times
\mathbb{C}^{m+d} \times \mathbb{C}^{rm} \rightarrow \C^{n+e}$$
are holomorphic on an open neighborhood of $$J^{2rk}(\mathbb{C}^{m+d},\mathbb{C}^{n+e})_{(0,0)}
\times J^{2rk}(\mathbb{C}^{m+d},\mathbb{C}^{n+e})_{(0,0)} \times
\mathbb{C}^{m+d} \times \mathbb{C}^{rm}.$$ Now
choose $x_0 \in \mathbb{C}^{rm}$ such that
$\hat{\delta}(t):=\delta(tx_0) \not\equiv 0$, for $t \in
\mathbb{C}$. As $H(Z)$ is independent of $x$, we can replace
$x=tx_0$ in (\ref{eq:eq57}). There exists a smallest integer $l$ such that
$\frac{d^l}{dt^l} \hat{\delta}(0) \neq 0$. To make our
calculations easier, consider a holomorphic change of variable
$\lambda=h(t)$ near the origin in $\mathbb{C}$, where $h$ is
determined by $\delta(tx_0)=\lambda^l$. So we now have
\begin{equation} \label{eq:eq49}
{\hat{\hat{\Xi}}}^{\mu,\nu}_{2r}\left(j_0^{2rk}H, j_0^{2rk}
\widetilde{H}, \frac{Z}{\lambda^l} ,\lambda \right) :=
\hat{\Xi}^{\mu,\nu}_{2r}\left(j_0^{2rk}H, j_0^{2rk}\widetilde{H}, \frac{Z}{\lambda^l} ,x_0h^{-1}(\lambda) \right) \equiv H(Z).
\end{equation}
Observe that the components of $\hat{\hat{\Xi}}^{\mu,\nu}_{2r}$
are in $$\mathcal{R}_0\Big(J^{2rk}(\mathbb{C}^{m+d},\mathbb{C}^{n+e})_{(0,0)} \times
 J^{2rk}(\mathbb{C}^{m+d},
 \mathbb{C}^{n+e})_{(0,0)}\times \mathbb{C}^{m+d} \times \mathbb{C} ,$$
   $$J^{2rk}(\mathbb{C}^{m+d},\mathbb{C}^{n+e})_{(0,0)} \times J^{2rk}(\mathbb{C}^{m+d},\mathbb{C}^{n+e})_{(0,0)} \Big) .$$
   To conclude the proof of the theorem, we expand the left hand side of (\ref{eq:eq49}) as a Laurent series in $\lambda$.
   Since $H(Z)$ is independent of $\lambda$, we can let $H(Z)$ be the constant term of the Laurent series. By Lemma
   \ref{thm:lemma2.11} and the form of $\Xi^{\mu,\nu}_{2r}$ given in (\ref{eq:niceForm}),
   we see that this is exactly of the form (\ref{eq:eq2.10}).

 Applying Lemma \ref{thm:lemma2.3}, we see that $\left(\bar{\tilde{H}},\bar{H}  \right)$
sends $\mathcal{M}$ into $\mathcal{M}'$ and satisfies condition
$D_{\nu\mu}$. From (\ref{eq:eq2.11}), we have
\begin{equation}
\bar{\tilde{H}}(Z) = \Phi^{\nu,\mu} \left(Z,
j_0^{K}\big(\bar{\tilde{H}} \big), j_0^K(\bar{H})   \right).
\end{equation}
Take the complex conjugate of this entire equation, and (\ref{eq:eq2.12}) follows.
%\hspace{1.222in}  $\Box$
$\Box$

\subsection{Reformulation of Theorem \ref{thm:Theorem1.4} }
 %$\,$ \vspace{.06in}

\begin{theorem2.11} \label{thm:Theorem2.12}
Let $M$ and $M'$ be as in Theorem \ref{thm:Finitejet}. Then there exists a positive
integer $L$, depending only on $M$ and $M'$, such that for each
$\alpha=(\alpha_1,\ldots,\alpha_n)$ with $1\leq \alpha_1 < \ldots <
\alpha_n \leq m$ and each $\beta=(\beta_1,\ldots,\beta_n)$ with $1
\leq \beta_1< \ldots<\beta_n \leq m$, there exist  $\C^{2n+2e}$-valued holomorphic
functions $\Phi_1^{\alpha,\beta}$ and $\Phi_2^{\alpha,\beta}$ defined on
an open subset of $\mathbb{C}^{m+d} \times \C^{m+d} \times
J^L(\mathbb{C}^{m+d},\mathbb{C}^{n+e})_{(0,0)}$ such that if
$\mathcal{H}$ is an HSPM satisfying condition $D_{\mu\nu}$, then
\begin{equation}  \label{eq:eq2.61}
\mathcal{H}(Z,\zeta)=\big(H(Z),\widetilde{H}(\zeta)\big) =
\Phi_1^{\mu,\nu}\big(Z,\zeta,j_0^L H\big) ,
\end{equation}
\begin{equation} \label{eq:eq2.62}
\mathcal{H}(Z,\zeta)=\big(H(Z),\widetilde{H}(\zeta)\big) =
\Phi_2^{\mu,\nu}\big(Z,\zeta,j_0^L \widetilde{H}\big) .
\end{equation}
\end{theorem2.11}
\begin{proof}  We will prove (\ref{eq:eq2.61}), and the proof
of (\ref{eq:eq2.62}) follows similarly. We will show inductively that there exist
$\C^{n+e}$-valued holomorphic functions $B_s^{\alpha,\beta,\gamma}$ defined on an open subset of
$J^{ks+|\gamma|}(\C^{m+d},\C^{n+e})_{(0,0)} \times \mathbb{C}^{ms}$, such that
\begin{equation} \label{eq:eq2.63}
\partial^\gamma H\big(v^s(t^0,\ldots,t^{s-1})\big) = B_s^{\mu,\nu,\gamma} \big(j_0^{ks+|\gamma|} \mathcal{G},t^0,\ldots,t^{s-1}\big) ,
\end{equation}
\begin{equation} \label{eq:eq2.64}
\partial^\gamma \widetilde{H}\big(\overline{v^s}(t^0,\ldots,t^{s-1})\big) =
\overline{B_s^{\nu,\mu,\gamma}} \big(j_0^{ks+|\gamma|}
\mathcal{G}',t^0,\ldots,t^{s-1}\big) ,
\end{equation}
where $\mathcal{G}=H$ and $\mathcal{G}'=\widetilde{H}$ if $s$ is even, and $\mathcal{G}=\widetilde{H}$ and $\mathcal{G}'={H}$ if $s$ is odd.

For $s=1$, we see that (\ref{eq:eq2.63}) and (\ref{eq:eq2.64}) hold true by letting
$(Z,\zeta)=\big( (z,0),0 \big)$ in (\ref{eq:eq15}) and $(Z,\zeta)=\big(0,(\chi,0)\big)$ in
(\ref{eq:eq16}). For some $s>1$, assume (\ref{eq:eq2.63}) and (\ref{eq:eq2.64}) hold for $s-1$. Assume, without loss of generality,
that $s$ is even (a similar proof works for $s$ odd).
 As $\big
(v^s(t^0,\ldots,v^{s-1}),
\overline{v^{s-1}}(t^1,\ldots,t^{s-1})\big) \in \mathcal{M}$, we see
from (\ref{eq:eq15}) that
\begin{equation} \label{eq:eq2.65}
\partial^\beta H\big(v^s(t^0,\ldots,t^{s-1})\big) \equiv \Psi_\beta^\nu\Big( v^s(t^0,\ldots,t^{s-1}),
\overline{v^{s-1}}(t^1,\ldots,t^{s-1}), j^{k+|\beta|}_{\overline{v^{s-1}}(t^1,\ldots,t^{s-1})} \tilde{H} \Big) .
\end{equation}
Using (\ref{eq:eq2.64}), we see then that
$$\partial^\beta H\big(v^s(t^0,\ldots,t^{s-1})\big) \equiv $$
\begin{equation}  \label{eq:eq2.65a}
\Psi_\beta^\nu \Big( v^s(t^0,\ldots,t^{s-1}),
\overline{v^{s-1}}(t^1,\ldots,t^{s-1}), \Big(
\overline{B_{s-1}^{\nu,\mu,\gamma}} \big(j_0^{k(s-1)+|\gamma|}
H,t^1,\ldots,t^{s-1}\big) \Big)_{|\gamma|\leq k+|\beta|} \Big).
\end{equation}
Now define $B_s^{\mu,\nu,\beta} \big(\Lambda,t^0,\ldots,t^{s-1}\big)$ to be the right hand side of (\ref{eq:eq2.65a}),
with the jets of $H$ replaced by the appropriate corresponding coordinates of $\Lambda$.

Using Lemma \ref{thm:lemma2.3}, we see that $\left(\bar{\tilde{H}},\bar{H}\right)$
satisfies condition $D_{\nu\mu}$ and sends $\mathcal{M}$ into
$\mathcal{M}'$. So, we have from (\ref{eq:eq2.63})
\begin{equation}
\partial^\gamma \bar{\tilde{H}}\big(v^s(t^0,\ldots,t^{s-1})\big) = B_s^{\nu,\mu,\gamma} \big(j_0^{ks+|\gamma|} \bar{\tilde{H}},t^0,\ldots,t^{s-1}\big).
\end{equation}
Taking the complex conjugate of both sides gives us (\ref{eq:eq2.64}).

Let $r$ be as given in Lemma \ref{thm:lemma2.9}. We know from (\ref{eq:eq2.63}) and (\ref{eq:eq2.64})
that
\begin{equation} \label{eq:eq2.68}
 H\big(v^{2r}(t^0,\ldots,t^{2r-1})\big) = B_{2r}^{\mu,\nu,0} \big(j_0^{2kr} H ,t^0,\ldots,t^{2r-1}\big) ,
\end{equation}
\begin{equation} \label{eq:eq2.69}
 \widetilde{H}\big(\overline{v^{2r+1}}(t^0,\ldots,t^{2r})\big) = \overline{B_{2r+1}^{\nu,\mu,0}} \big(j_0^{2kr+k} H ,t^0,\ldots,t^{2r}\big) .
\end{equation}
As $v^{l+1}(t^0,\ldots,t^{l-1},0) = v^l(t^0,\ldots,t^{l-1})$ for any
positive integer $l$, we see from Lemma \ref{thm:lemma2.9} that the matrix
\begin{equation}
 \left( \frac{\partial v^{2r+1} }{\partial
(t^0,t^{r+1},t^{r+2},\ldots,t^{2r-1})}
(0,x^1,\ldots,x^{r-1},x^r,x^{r-1},\ldots,x^1,0) \right)
\end{equation}
has rank $m+d$ for all $(x_1, \ldots, x_r) \in U \backslash V$, for $U \subseteq \C^{rm}$ an open neighborhood of the origin and $V$ a proper
holomorphic subvariety of $U$, and we also see that
\begin{equation}
v^{2r+1} (0,x^1,\ldots,x^{r-1},x^r,x^{r-1},\ldots,x^1,0) \equiv 0 .
\end{equation}
We can now use (\ref{eq:eq2.68}) and (\ref{eq:eq2.69}) to obtain (\ref{eq:eq2.61}) and (\ref{eq:eq2.62}) by following
exactly the proof of Theorem \ref{thm:Finitejet}.
\end{proof}

\subsection{Reformulation of Theorem \ref{thm:algebraic} }
 %$\,$ \vspace{.06in}

\begin{algebraicII} \label{thm:algebraicII}
Let $M$ and $M'$ be as in Theorem \ref{thm:Finitejet}, and assume that $M$ and $M'$ are real algebraic. If
$\mathcal{H}$ is an HSPM satisfying condition $D_{\mu\nu}$
for some $\mu$ and $\nu$, then
$\mathcal{H}$ is holomorphic algebraic.
\end{algebraicII}

\begin{proof}
An inspection of the proof of Lemma \ref{thm:lemma2.4} shows that the $\phi_\beta^\alpha$ as given in (\ref{eq:eq4})
are holomorphic algebraic (as $M$ is real  algebraic).  When solving the system of equations in (\ref{eq2.40z}), apply
the \emph{algebraic} implicit function theorem  to see that $B^\nu$ as given in (\ref{eq:eq18})
 and (\ref{eq:eq22}) is holomorphic algebraic (as $M'$
is real algebraic).  Thus, an inspection of the proof of Lemma \ref{thm:lemma2.5} shows that the $\Psi^\alpha_\beta$ as
given in (\ref{eq:eq14}) are holomorphic algebraic.  An examination  of the proof of Lemma \ref{lemmalemma} then reveals that
the $\Xi_s^{\alpha,\beta}$ as given in (\ref{eq:niceForm}) are holomorphic algebraic.  Finally, in the proof
of Theorem \ref{thm:Finitejet}, choose $x_0$  sufficiently small and satisfying $\delta(x_0) \neq 0$, and substitute $x=x_0$ in (\ref{eq:eq56}).
By Lemma \ref{thm:lemma2.10}, we see then that $H(Z)$ is holomorphic algebraic.
Similarly, $\widetilde{H}(\zeta)$ is holomorphic algebraic.
\end{proof}

\section{Proofs of Main Results}    \label{section:2.4}
In Section 1, we presented Theorem {\ref{thm:Theorem1.1}, Corollary \ref{thm:cor1.2}, and
Corollary \ref{thm:cor1.3}, all of which follow naturally from Theorem \ref{thm:Finitejet}.
We also presented Theorem \ref{thm:Theorem1.4}, which is a direct result of Theorem \ref{thm:Theorem2.12},
and Theorem \ref{thm:algebraic}, which is a direct result of Theorem \ref{thm:algebraicII}.
In this section, we provide their proofs. First we make the following observations.

\begin{segSubObs}  \label{segSub}
If $M \subseteq \C^N$ and $M' \subseteq \C^{N'}$ are submanifolds of codimensions $d$ and $d'$, respectively, given in normal coordinates by $w=Q(z,\X,\T)$ and $w'=Q'(z',\X',\T')$, respectively, then an HSPM $\HH=(f,g,\tilde{f},\tilde{g})$ sending $(\M,0)$ into $(\M',0)$ is Segre submersive at $0$ if and only if
the matrices $\big(f_z(0)\big)$ and $\big(\f_\X(0)\big)$ have rank $N'-d'$.
 This follows from the fact that a basis for the antiholomorphic vectors tangent to $M$ (resp., $M'$) at $0$ is given by $\big\{
\frac{\partial}{\partial \bar{z}_j} : 1 \leq j \leq N-d\big\}$ $\big($resp., $\big\{
\frac{\partial}{\partial \bar{z}'_j} : 1 \leq j \leq N'-d'\big\}\big)$,
and a basis for the holomorphic vectors tangent to $M$ (resp., $M'$) at $0$ is given by $\big\{
\frac{\partial}{\partial {z}_j} : 1 \leq j \leq N-d \big\}$ $\big($resp., $\big\{
\frac{\partial}{\partial {z}'_j} : 1 \leq j \leq N'-d'\big\}\big)$,
coupled with the fact that $g_{z_j}(0) =\g_{\X_j}(0)= 0$ for $j=1,\ldots,N-d$.
\end{segSubObs}

\begin{nextObs} \label{nextObs}
For $p \in \C^N$, let $\phi : (\C^N,0) \rightarrow (\C^N,p)$ be a biholomorphism near $0$, and
for $p' \in \C^{N'}$, let $\phi' : (\C^{N'},0) \rightarrow (\C^{N'},p')$ be a biholomorphism near $0$.
Then for any nonnegative $l$, there exist vector-valued polynomial functions $F_l$ and $G_l$ such that if
$h:(\C^N,0) \rightarrow (\C^{N'},0)$ is any holomorphic map, and $\tilde{h} : (\C^N,p) \rightarrow (\C^{N'},p')$
 is given by $\tilde{h} := \phi' \circ h \circ \phi^{-1}$, then
$j_p^l \tilde{h} = F_l(j_0^l h)$ and $j_0^l h = G_l(j_p^l \tilde{h})$.
\end{nextObs}

\subsection{Proof of Theorem \ref{thm:Theorem1.1}}
 %$\,$ \vspace{.06in}

Theorem \ref{thm:Theorem1.1} follows  from Theorem \ref{thm:Finitejet}, Observation \ref{segSub}, and Observation \ref{nextObs}.
We leave the details to the reader. $\Box$
%%\hspace{4.85in} $\Box$

 \subsection{Proof of Corollary \ref{thm:cor1.2}}
  %$\,$ \vspace{.06in}

 Without loss of generality, assume $p=0$.
 As $M=M'$, it is clear from the statement of Theorem \ref{thm:Finitejet} that we can can
 choose $r=1$ in Theorem \ref{thm:Theorem1.1}. Do so, and define $\Phi := \Phi_1$ as given in (\ref{phi1}).
 %(The fact that $r=1$ in Theorem \ref{thm:Theorem1.1} is clear from the statement of Theorem \ref{thm:Finitejet}.)
It then follows from Theorem \ref{thm:Theorem1.1}
that $\eta_0^K$ is continuous and injective on $\Aut$.  To show that $(\eta_0^K)^{-1}$
is continuous on $\eta_0^K(\Aut)$, let $\Lambda_j,\widetilde{\Lambda}_j,\Lambda_0,\widetilde{\Lambda}_0 \in G_0^k(\C^{N})$ and assume
$(\Lambda_j,\widetilde{\Lambda}_j) \in \eta_0^K(\Aut)$ converges to
$(\Lambda_0,\widetilde{\Lambda}_0) \in \eta_0^K(\Aut)$.
 Theorem \ref{thm:Theorem1.1} tells us that
$(\eta_0^K)^{-1}(\Lambda_j,\widetilde{\Lambda}_j) =\big( \Phi(Z,\Lambda_j,\widetilde{\Lambda}_j) , \overline{ \Phi}(\zeta,\widetilde{\Lambda}_j,{\Lambda}_j) \big)$
and $(\eta_0^K)^{-1}(\Lambda_0,\widetilde{\Lambda}_0) = \big( \Phi(Z,\Lambda_0,\widetilde{\Lambda}_0) , \overline{ \Phi}(\zeta,\widetilde{\Lambda}_0,{\Lambda}_0) \big)$.
 Furthermore, $\Phi(Z,\Lambda, \widetilde{\Lambda})$ is holomorphic in a neighborhood of $(0,\Lambda_0,\widetilde{\Lambda}_0)$ and thus continuous,
 and $\overline{\Phi}(\zeta,\widetilde{\Lambda}, {\Lambda})$ is holomorphic in a neighborhood of $(0,\widetilde{\Lambda}_0,{\Lambda}_0)$ and thus continuous.
 %Thus, $\big( \Phi(Z,\Lambda_j,\widetilde{\Lambda}_j) , \overline{ \Phi}(\zeta,\widetilde{\Lambda}_j,{\Lambda}_j) \big)$
 %converges uniformly to $\big( \Phi(Z,\Lambda_0,\widetilde{\Lambda}_0) , \overline{ \Phi}(\zeta,\widetilde{\Lambda}_0,{\Lambda}_0) \big)$
 %on compact neighborhoods of 0.
 Therefore, $(\eta_0^K)^{-1}(\Lambda_j,\widetilde{\Lambda}_j)$ converges to $(\eta_0^K)^{-1}(\Lambda_0,\widetilde{\Lambda}_0)$.
It follows then that $\eta_0^K$
 is a homeomorphism from $\Aut$ onto $\eta_0^K(\Aut)$.

We now show that $\eta_0^K \big(
\text{Aut}_{\mathbb{C}}(\mathcal{M},0) \big)$ is a closed,
holomorphic algebraic submanifold of $G_0^K(\mathbb{C}^{N})
\times G_0^K(\mathbb{C}^{N})$.
Let $\rho(Z,\bar{Z})$ be a defining function for $M$. Write $\zeta = ({\zeta}_1, {{\zeta}}_2) \in \C^{N-d}
\times \C^d$, where $d$ is the codimension of $M$.  After a possible rearrangement of coordinates, as $M$ is generic, there exists a
holomorphic map $\theta: \C^N \times \C^{N-d} \rightarrow \C^d$ satisfying $\theta(0)=0$ such that for all $Z$
and ${\zeta}_1$ sufficiently close to 0,
$(Z,{\zeta}_1, \theta(Z,{\zeta}_1)) \in \mathcal{M}$.
Given $(\Lambda_0,\widetilde{\Lambda}_0) \in G_0^K(\mathbb{C}^{N}) \times
G_0^K(\mathbb{C}^{N})$,
 $(\Lambda_0,\widetilde{\Lambda}_0) \in \eta_0^K \big( \text{Aut}_{\mathbb{C}}(\mathcal{M},0) \big)$ if and only if
 the following three conditions hold:
\begin{equation} \label{cond1}
\Lambda_0 = \big( \, S_\gamma(\Lambda_0,\widetilde{\Lambda}_0) \, \big)_{|\gamma| \leq K}
\end{equation}
\begin{equation} \label{cond2}
\tilde{\Lambda}_0 = \big(\,\, \overline{S_\gamma}(\widetilde{\Lambda}_0,{\Lambda}_0) \, \big)_{|\gamma| \leq K}
\end{equation}
\begin{equation} \label{cond3}
\rho\left( \Phi(Z,\Lambda_0,\widetilde{\Lambda}_0) , \overline{\Phi}\big({\zeta}_1, \theta(Z, \zeta_1),\widetilde{\Lambda}_0, \Lambda_0\big)\right) = 0,
\end{equation}
where $S_\gamma$
are the rational coefficients in the Taylor expansion given in (\ref{phi1}).
 Equations (\ref{cond1}) and (\ref{cond2}) can be expressed as a finite set of polynomial
equations in $\Lambda_0$ and $\widetilde{\Lambda}_0$ as each $S_\gamma$ is rational. Equation (\ref{cond3}) can be expressed as an infinite
set of polynomial equations in $\Lambda_0$ and $\widetilde{\Lambda}_0$.  This can be seen by noting that $\Phi(0,\Gamma,\widetilde{\Lambda}) \equiv 0$ and $\theta(0)=0$,
and by noting the form of $\Phi$ given in Theorem \ref{thm:Theorem1.1}.

 Thus, we see that $\eta_0^K \big(
\text{Aut}_{\mathbb{C}}(\mathcal{M},0) \big)$ is a closed, holomorphic algebraic
subvariety of the space
%$\,\,\,\,\,\,\,$
$G_0^K(\mathbb{C}^{N}) \times G_0^K(\mathbb{C}^{N})$
as it is given by the vanishing of a set of polynomial equations. To
see that it is actually a submanifold, we first note that it is a
subgroup of  $G_0^K(\mathbb{C}^{N}) \times
G_0^K(\mathbb{C}^{N})$ as multiplication can be defined in the
following way: given any $(\Lambda_1,\widetilde{\Lambda}_1),$
$(\Lambda_2,\widetilde{\Lambda}_2) \in \eta_0^K \big(
\text{Aut}_{\mathbb{C}}(\mathcal{M},0) \big)$, let $\mathcal{H}_1$
and $\mathcal{H}_2$, respectively, be the corresponding
automorphisms in $\text{Aut}_{\mathbb{C}}(\mathcal{M},0)$. Now
compose $\mathcal{H}_1$ and $\mathcal{H}_2$, and apply $\eta_0^K$ to
this composition. Under this multiplication, $\eta_0^K \big(
\text{Aut}_{\mathbb{C}}(\mathcal{M},0) \big)$ is a closed subgroup
of the Lie group
$G_0^K(\mathbb{C}^{N}) \times G_0^K(\mathbb{C}^{N})$, and is thus a Lie subgroup (see, for example, \cite{V74}).
%%\hspace{0.6in}
$\Box$

\subsection{Proof of Corollary \ref{thm:cor1.3}}
 %$\,$ \vspace{.06in}

Before proving Corollary \ref{thm:cor1.3}, we present a simple lemma which
involves only basic linear algebra.
\begin{lemma3.1} \label{thm:lemma2.4.1}
Let $A=(a_{ij})$ be a $d \times d$ invertible matrix, where
each $a_{ij} \in \C$. Let $b_1,\ldots,b_d \in \C$. Let $B_1$ be the matrix gotten
by replacing row $m$ of $A$ with $(a_{m1}+b_1, \ldots, a_{md}+b_d)$
and $B_2$ be the matrix gotten by replacing row $m$ of $A$ with
$(a_{m1}-b_1, \ldots, a_{md}-b_d)$. Then at least one of $B_1$ or
$B_2$ is invertible.
\end{lemma3.1}
\begin{proof} Without loss of generality, assume $m=1$.
Let $A_n := (-1)^{n+1}\det M_n$, where $M_n$ is
the $(d-1) \times (d-1)$ matrix gotten by deleting the first row and
$n^{th}$ column of $A$. Assume that $\det B_1 = \det B_2 =0$. Then
expanding along the first row of $B_1$  gives
\begin{equation}  \label{Aa123}
(a_{11} + b_1)A_1 + \ldots + (a_{1d}+b_d)A_d = 0,
\end{equation}
and expanding along the first row of $B_2$ gives
\begin{equation} \label{Aa234}
(a_{11} - b_1)A_1 + \ldots + (a_{1d}-b_d)A_d = 0 .
\end{equation}
Adding (\ref{Aa123}) and (\ref{Aa234})  gives $2a_{11}A_1+\ldots+2a_{1d}A_d=0$. However, this implies that $\det A=0$, a contradiction.
\end{proof}

We now prove Corollary \ref{thm:cor1.3}.
 Without loss of generality, assume $p=0$. Let $r(\Lambda,\bar{\Lambda}) = \big( r_1(\Lambda,\bar{\Lambda}), \ldots, r_s(\Lambda,\bar{\Lambda}) \big)$
 be a defining function for $j_0^K\big(\text{Aut}(M,0)\big)$ as a real algebraic
submanifold of $G_0^K(\mathbb{C}^{N})$,
where $\Lambda \in G_0^K(\C^N)$
 (we refer the reader to Remark \ref{remark2.2.4}).
 The complexification of this submanifold, $\C \big\{ j_0^K\big(\text{Aut}$ $(M,0)\big)\big\}$, is thus a complex submanifold of
  $G_0^K(\mathbb{C}^{N}) \times G_0^K(\mathbb{C}^N)$, given by the vanishing of $r(\Lambda,\widetilde{\Lambda})$,
  where $\widetilde{\Lambda} \in G_0^K(\C^N)$. Let $\rho(Z,\bar{Z})$ be a
   defining function for $M$.  As $M=M'$, it is clear from the statement of Theorem \ref{thm:Finitejet} that we can can
 choose $r=1$ in Theorem \ref{thm:Theorem1.1}. Do so, and define $\Phi := \Phi_1$ as given in (\ref{phi1}).
From Theorem \ref{thm:Theorem1.1}, we see that for any $\Lambda \in
G_0^K(\mathbb{C}^{N})$,
\begin{equation}
\rho\big( \Phi(Z,\Lambda,\bar{\Lambda}) ,
\overline{\Phi}(\bar{Z},\bar{\Lambda},{\Lambda}) \big) = A(Z,\Lambda,
\bar{Z},\bar{\Lambda})r(\Lambda,\bar{\Lambda}) +
B(Z,\Lambda,\bar{Z},\bar{\Lambda})\rho(Z,\bar{Z}),
\end{equation}
where $A$ is a real analytic $d \times s$ matrix, and $B$ is a real analytic
$d \times d$ matrix. Complexify to get:
\begin{equation} \label{eq:eq2.73}
\rho\big( \Phi(Z,\Lambda,\widetilde{\Lambda}) ,
\overline{\Phi}(\zeta,\widetilde{\Lambda},{\Lambda}) \big) = A(Z,\Lambda,
\zeta,\widetilde{\Lambda})r(\Lambda,\widetilde{\Lambda}) +
B(Z,\Lambda,\zeta,\widetilde{\Lambda})\rho(Z,\zeta) .
\end{equation}
Notice that (\ref{eq:eq2.73}) gives us exactly what we want. This equation says that if
$(\Lambda,\widetilde{\Lambda}) \in \C \big\{
j_0^K\big(\text{Aut}(M,0)\big)\big\}$, then $
\big(\Phi(Z,\Lambda,\widetilde{\Lambda})  ,
\overline{\Phi}(\zeta,\widetilde{\Lambda},{\Lambda})  \big)   \in
\text{Aut}_{\mathbb{C}}(\mathcal{M},0)$. Now we need only show that
\begin{equation}
\eta_0^K \big(\Phi(Z,\Lambda,\widetilde{\Lambda})  ,
\overline{\Phi}(\zeta,\widetilde{\Lambda},{\Lambda})  \big)  =
(\Lambda,\widetilde{\Lambda}).
\end{equation}
 We  have the equations:
$$   \Big( \partial^\alpha_Z \Phi(0,\Lambda,\bar{\Lambda}) \Big)_{ |\alpha| \leq K} =
\Lambda +  C(\Lambda,\bar{\Lambda})r(\Lambda,\bar{\Lambda}) ,$$
\begin{equation}
 \Big( \partial^\alpha_{\bar{Z}} \overline{\Phi}(0,\bar{\Lambda},{\Lambda})  \Big)_{ |\alpha| \leq K} =
 \bar{\Lambda} + \overline{C}(\bar{\Lambda},{\Lambda})r(\Lambda,\bar{\Lambda}),
\end{equation}
for $C$ a real analytic matrix. Complexify these to get:
$$   \Big( \partial^\alpha_Z \Phi(0,\Lambda,\widetilde{\Lambda}) \Big)_{|\alpha| \leq K} =
\Lambda + C(\Lambda,\widetilde{\Lambda})r(\Lambda,\widetilde{\Lambda}) ,$$
\begin{equation}
 \Big( \partial^\alpha_\zeta \overline{\Phi}(0,\widetilde{\Lambda},{\Lambda})  \Big)_{ |\alpha| \leq K} =
 \widetilde{\Lambda} +  \overline{C}(\widetilde{\Lambda},{\Lambda})r(\Lambda,\widetilde{\Lambda}),
\end{equation}
and the first part of Corollary \ref{thm:cor1.3} is proved.

As we are assuming $p=0$, we take $Id=Id'$ in Corollary \ref{thm:cor1.3}.
To prove the second part of Corollary \ref{thm:cor1.3},
 first we show that near $(Id,Id)$,
$\eta^K_0 \big(\Aut\big)$ is a complexified submanifold. In other words,
$\eta^K_0 \big(\Aut\big)$
$=\mathbb{C} R$, where $R$ is a real
submanifold of $G_0^K(\C^{N})$ (here, $\C R$ denotes the
complexification of $R$). We know from Corollary \ref{thm:cor1.2} that
$\eta_0^K(\Aut)$ is a complex submanifold of $G_0^{K}(\C^{N})
\times G_0^{K}(\C^{N})$.  Near $(Id,Id)$, let
$\hat{s}_1(\Lambda,\widetilde{\Lambda}),\ldots,\hat{s}_t(\Lambda,\widetilde{\Lambda})$
be  defining functions for $\eta_0^K(\Aut)$. We will assume without loss of generality that these
functions are defined on a ball $B$ of sufficiently small radius centered at $(Id, Id)$; this way if
$(\Gamma,\Lambda)$ is a point in $B$, then so is
$(\Lambda,\Gamma)$ and $(\bar{\Lambda},\overline{\Gamma})$. Now we set
\begin{equation} \label{q111}
s_j(\Lambda,\widetilde{\Lambda}) := \hat{s}_j(\Lambda,\widetilde{\Lambda})+\bar{\hat{s}}_j(\widetilde{\Lambda},{\Lambda})
\end{equation}
 $$\text{\emph{or}} $$
\begin{equation}  \label{q222}
s_j(\Lambda,\widetilde{\Lambda}) :=
i\hat{s}_j(\Lambda,\widetilde{\Lambda})-i\bar{\hat{s}}_j(\widetilde{\Lambda},{\Lambda}).
\end{equation}

 We choose between options (\ref{q111}) and (\ref{q222}) as follows. Start with $j=1$. From Lemma \ref{thm:lemma2.4.1},
 we can replace $\hat{s}_1$ with one of the above $s_1$,
 and in at least one case the differentials of $s_1,\hat{s}_2,\ldots,\hat{s}_t$ will be linearly independent near $(Id,Id)$. Choose $s_1$ so that this
 is the case. Now do the same thing for $j=2$, then $j=3$, and so forth. Let $\mathcal{R}$ be the submanifold defined by
 $s_1(\Lambda,\tilde{\Lambda})=\ldots=s_t(\Lambda,\tilde{\Lambda})=0$. If $(\Lambda,\tilde{\Lambda}) \in \eta_0^K(\Aut)$,
 then from Lemma \ref{thm:lemma2.3}, $(\bar{\tilde{\Lambda}},\bar{\Lambda}) \in \eta_0^K(\Aut)$. Thus,
 $\hat{s}_j(\bar{\tilde{\Lambda}},\bar{\Lambda})=0$,
  which  implies that
  $\bar{\hat{s}}_j (\tilde{\Lambda},\Lambda)=0$. In other words, near $(Id,Id)$, $\eta_0^K(\Aut) \subseteq \mathcal{R}$.
  But these two submanifolds have equal dimensions. So we see that, in fact, $\eta_0^K(\Aut) = \mathcal{R}$ near $(Id, Id)$.

Now we need only show that $\mathcal{R} = \C R$ for some
real submanifold $R \subseteq G_0^{K}(\C^{N})$, and we will have proved
our claim. Let
\begin{equation}
R := \{\Lambda: s_1(\Lambda,\bar{\Lambda})=\ldots=s_t(\Lambda,\bar{\Lambda})=0 \} .
\end{equation}
Clearly $R$ is a nonempty set as it contains the point $\Lambda = Id$.
 As each $s_j$ is a real function, and the differentials
 of $s_1,\ldots,s_t$ are linearly independent,
it follows that $R$ is a real submanifold.

From Theorem \ref{thm:Theorem1.1}, we see that if $(H,\widetilde{H}) \in \Aut$
and  $j_0^K(\widetilde{H}) = j_0^K(\, \overline{H})$,  we must have
$\widetilde{H}=\overline{H}$. Thus, near $(Id,Id)$:
\begin{equation}
\C \big\{j_0^K(\text{aut}(M,0))\big\} \cap
\{\widetilde{\Lambda}=\bar{\Lambda} \} = \eta_0^K({\Aut}) \cap
\{\widetilde{\Lambda}=\bar{\Lambda} \} = \C R \cap
\{\widetilde{\Lambda}=\bar{\Lambda} \} ,
\end{equation}
implying that $j_0^K(\text{aut}(M,0))=R$. Thus their complexifications must be equal as well.
 That is, near $(Id,Id)$, $\C \big\{j_0^K(\text{aut}(M,0))\big\}  = \eta_0^K(\Aut)$.
 But both of these are algebraic holomorphic submanifolds.
 So if they are equal near $(Id,Id)$, then using the notation given in the statement of this corollary, we must have $\BB = \CC$.

The third part of the corollary comes from the fact that $\eta_0^K(\Aut)$ is a Lie subgroup.
  Thus, each of its connected components is a coset of $\BB$.
  Since $\C \big\{ j_0^K\big(\text{Aut}(M,0)\big)\big\} \subseteq \eta_0^K(\Aut)$, and they are both algebraic holomorphic submanifolds,  each component of $\C \big\{ j_0^K\big(\text{Aut}(M,0)\big)\big\} $ is exactly equal to one of the components of  $\eta_0^K(\Aut)$.
  Algebraicity implies that there are finitely many such components.
\hspace{1.6in}    $\Box$

\subsection{Proof of Theorem \ref{thm:Theorem1.4}}

Theorem \ref{thm:Theorem1.4} follows  from Theorem \ref{thm:Theorem2.12}, Observation \ref{segSub}, and Observation \ref{nextObs}.
We leave the details to the reader.
%\hspace{4.89 in}     $\Box$
$\Box$

\subsection{Proof of Theorem \ref{thm:algebraic}}

As $M$ and $M'$ are real algebraic, they have real analytic algebraic defining functions .
%$\rho(Z,\bar{Z})$ and $\rho'(Z',\bar{Z}')$, respectively.
When $M$ and $M'$ are expressed in normal coordinates, the new defining
functions can also be chosen to be real analytic algebraic. This follows by using the
algebraic implicit function theorem in the derivation of the the new defining functions (for
precise details on deriving normal coordinates and the algebraic implicit function theorem, see \cite{BER99a}).
Furthermore, if $\tilde{Z} = \varphi(Z)$ is a holomorphic algebraic change of coordinates, then $\varphi^{-1}$ is
a holomorphic algebraic function (this is also a direct consequence of the algebraic implicit function theorem).
  Thus, Theorem \ref{thm:algebraic} now
follows  from Theorem \ref{thm:algebraicII} and Observation \ref{segSub}.
\hspace{4.89 in}     $\Box$

\section{Examples: HSPMs and Automorphism Groups}    \label{section:Examples}

%PUT AN INTRODUCTION (IN FIRST SECTION, ETC.)

For $n > 1$, there exist $M, M' \subseteq \C^{n+1}$
defined near 0 such that there exist no holomorphic
maps $H$ satisfying:
\begin{equation} \label{eq:realmap}
H \text{ is invertible near 0}, \,\, H(M) \subseteq M', \,\, H(0)=0,
\end{equation}
 yet there exist  HSPMs  satisfying:
 \begin{equation} \label{eq:Segremap}
 \HH \text{ is invertible near 0}, \,\, \HH(\M) \subseteq \M', \,\, \HH(0)=0.
 \end{equation}

\begin{example31} \label{ex31}
For $n>1$, let $(z_1,\ldots,z_n,w)$ and $(z_1',\ldots,z_n',w')$ be coordinates on $\C^{n+1}$ and define
\begin{equation*}
M=\left\{\text{Im } w= \sum_{j=1}^n \epsilon_j |z_j|^2   \right\}  ,
\end{equation*}
 \begin{equation*}
 M'=\left\{\text{Im } w' = \sum_{j=1}^n \sigma_j |z_j'|^2   \right\} ,
 \end{equation*}
 where $\epsilon_j, \sigma_j \in \{-1,1\}$.
 Both $M$ and $M'$ are of finite type and finitely nondegenerate at 0.
  If $\left| \sum_j \epsilon_j \right| \neq  \left| \sum_j \sigma_j \right|$, then there are no holomorphic maps satisfying
 criteria (\ref{eq:realmap}).
 %(indeed, if $M$ and $M'$ have different Levi signatures, then they are not biholomorphically equivalent).
 %\emph{\textbf{WHY IS THIS TRUE?????????????????????????????????}} \\
 %\emph{\textbf{AT LEAST I THINK IT'S TRUE.......................}} \\
(Indeed,
$M$ and $M'$ have different Levi signatures at 0.)
% see, for example, \cite{Krantz}.)
 However, for $a, c_j \in \C\backslash\{0\}$, the family of maps given by
 \begin{equation*}
  \mathcal{H}(z,w,\X,\T)=
  \end{equation*}
  \begin{equation*}
  \left(\epsilon_1c_1z_1\,,\,\ldots,\epsilon_{n-1}c_{n-1}z_{n-1}\,,\,\epsilon_nc_nz_n\,,\,aw\,,\,
  \frac{a\sigma_1}{c_1}\chi_1\,,\,\ldots\,,\,
\frac{a\sigma_{n-1}}{c_{n-1}}\chi_{n-1}\,,\, \frac{a\sigma_n}{c_n}\chi_n\,,\,a\tau \right)
  \end{equation*}
   satisfy criteria (\ref{eq:Segremap}).
  \end{example31}

This can also occur in $\mathbb{C}^2$ as the next example illustrates.

\begin{example32}  \label{ex32}
Let $M,M' \subseteq \C^2$ be given by
\begin{equation*}
  M= \Big\{ \text{Im } w=|z|^2+ 2\text{Re}\big[z^4\bar{z}^2(1+i\text{Re }w)\big] \Big\},
\end{equation*}
\begin{equation*}
M'= \Big\{ \text{Im } w'=|z'|^2+ 2\text{Re}\big[z'^4\bar{z}'^2(1-i\text{Re }w')\big] \Big\}.
\end{equation*}
 Notice that $M$ and $M'$ are of finite type and finitely nondegenerate at 0.
It can be shown (\cite{CM}) that there are no maps $H$ satisfying criteria (\ref{eq:realmap}).
(Indeed, as $M$ and $M'$ are in Chern-Moser normal form, the fact that the coefficients $i$ and $-i$ are unequal
implies that there does not exist a holomorphic map $H$ satisfying criteria (\ref{eq:realmap}).)
 However, it easy to check that the HSPM
$\mathcal{H}(z,w,\X,\T)=(iz,-w,i\chi,-\tau)$ satisfies (\ref{eq:Segremap}).
\end{example32}

Now we will look at some examples of automorphism groups.  In Example \ref{ex4.1bb}, we find that
$\C \big\{j_0^K\big(\text{Aut}(M,0)\big)\big\}$  and $\eta_0^K(\Aut)$ are equal.
\begin{ex4.1} \label{ex4.1bb}
 Let $M$ be the \emph{Lewy hypersurface} of $\C^2$. It is given by
\begin{equation*}
M = \{\text{Im }w=|z|^2\} .
 \end{equation*}
 We note that $M$ is finitely nondegenerate and of finite type at 0.
 It can be shown  (see \cite{Angle} for the calculations) that every holomorphic
  Segre preserving automorphism of $\mathcal{M}$ at 0 is of the form
$$\HH(z,w,\X,\T) = \Bigg(  \frac{\alpha (z +\beta w)}{1-(\gamma + i\beta\tilde{\beta}) w -2i\tilde{\beta}z},
 \frac{\alpha \tilde{\alpha} w}{1-(\gamma + i\beta\tilde{\beta}) w -2i\tilde{\beta}z}, $$
\begin{equation} \label{eq2.78}
\frac{\tilde{\alpha}(\chi + \tilde{\beta}\tau)}{1-(\gamma -
i\beta\tilde{\beta}) \tau +2i{\beta}\chi},
\frac{\alpha\tilde{\alpha}\tau}{1-(\gamma - i\beta\tilde{\beta})
\tau +2i{\beta}\chi} \Bigg),
\end{equation}
where $\gamma,\beta,\tilde{\beta} \in
\mathbb{C}$ and $\alpha,\tilde{\alpha} \in \C \backslash \{0\}$. Also, every  automorphism of $M$ at 0 is of the form
\begin{equation} \label{eq2.79}
H(z,w) = \Bigg(  \frac{\alpha (z +\beta w)}{1-(\gamma + i|\beta|^2) w -2i\bar{\beta}z},   \frac{|\alpha|^2 w}{1-(\gamma + i|\beta|^2) w -2i\bar{\beta}z} \Bigg),
\end{equation}
%\end{equation}
where $\alpha \in \C \backslash \{0\}$, $\beta \in \mathbb{C}$, and $\gamma \in \mathbb{R}$.
The automorphisms in (\ref{eq2.79}) follow directly from
the automorphisms in (\ref{eq2.78}), but those in (\ref{eq2.79}) have actually been known for some time (see \cite{CM}).

We see from (\ref{eq2.78}) and (\ref{eq2.79}) that  $ \C \big\{
j_0^K\big(\text{Aut}(M,0)\big)\big\} = \eta_0^K(\Aut)$.
Indeed, let $(\Lambda^f_z, \ldots, \Lambda^f_{ww},\Lambda^g_z, \ldots, \Lambda^g_{ww},
\widetilde{\Lambda}^{\f}_{\X} , \ldots, \widetilde{\Lambda}^{\f}_{\T\T},
\widetilde{\Lambda}^{\g}_{\X} , \ldots, \widetilde{\Lambda}^{\g}_{\T\T})$ be
coordinates on $G_0^2(\C^2) \times G_0^2(\C^2)$, where
$\Lambda^f_{z^r w^s}$ corresponds to $\frac{\partial^{r+s}f}{\partial z^r \partial w^s}$,
$\Lambda^g_{z^r w^s}$ corresponds to $\frac{\partial^{r+s}g}{\partial z^r \partial w^s}$,
$\Lambda^{\f}_{\X^r \T^s}$ corresponds to $\frac{\partial^{r+s}\f}{\partial \X^r \partial \T^s}$, and
$\Lambda^{\g}_{\X^r \T^s}$ corresponds to $\frac{\partial^{r+s}\g}{\partial \X^r \partial \T^s}$.
Then (\ref{eq2.79}) implies that
$ \C \big\{j_0^2\big(\text{Aut}(M,0)\big)\big\}$ is given by
\begin{equation*}
\Bigg\{
\Lambda^g_w=\widetilde{\Lambda}^{\g}_\T=\Lambda^f_z\widetilde{\Lambda}^{\f}_\X, \Lambda^g_{ww}-\widetilde{\Lambda}^{\g}_{\T\T} =2i\Lambda^f_w\widetilde{\Lambda}^{\f}_\T,
 \Lambda^g_{zw} =2i\Lambda^f_z\widetilde{\Lambda}^{\f}_\T,
\Lambda^f_{zw}=\frac{\Lambda^g_{ww}}{\widetilde{\Lambda}^{\f}_\X}, \Lambda^f_{zz}=2i\frac{\Lambda^f_z\widetilde{\Lambda}^{\f}_\T}{\widetilde{\Lambda}^{\f}_\X},
\end{equation*}
\begin{equation*}
\Lambda^f_{ww}=\frac{\Lambda^g_{ww}\Lambda^f_w}{\Lambda^f_z\widetilde{\Lambda}^{\f}_\X} ,
 \widetilde{\Lambda}^{\g}_{\X\T}=-2i\widetilde{\Lambda}^{\f}_\X{\Lambda}^f_w,
\widetilde{\Lambda}^{\f}_{\X\T}=\frac{\widetilde{\Lambda}^{\g}_{\T\T}}{{\Lambda}^f_z}, \widetilde{\Lambda}^{\f}_{\X\X}=-2i\frac{\widetilde{\Lambda}^{\f}_\X{\Lambda}^f_w}{{\Lambda}^f_z},
\widetilde{\Lambda}^{\f}_{\T\T}=\frac{\widetilde{\Lambda}^{\g}_{\T\T}\widetilde{\Lambda}^{\f}_\T}{\widetilde{\Lambda}^{\f}_\X{\Lambda}^f_z} ,
\end{equation*}
\begin{equation}    \label{ssaass}
 \Lambda^g_z=\Lambda^g_{zz}=\widetilde{\Lambda}^{\g}_\X=\widetilde{\Lambda}^{\g}_{\X\X}=0 \Bigg\} .
\end{equation}
It follows from (\ref{eq2.78}) that $\eta_0^2(\Aut)$ is also given by (\ref{ssaass}).

What is more interesting, however, are submanifolds for which
$\C \big\{j_0^K\big(\text{Aut}(M,0)\big)\big\} \neq \eta_0^K(\Aut)$.
\end{ex4.1}

\begin{ex4.2} \label{ex4.2bb}
Let $M \subseteq \C^2$ be given by
\begin{equation*}
 M = \Big\{\text{Im } w =|z|^2+(\text{Re
} z^2)|z|^2 \Big\} .
\end{equation*}
Notice that $M$ is finitely nondegenerate and of finite type at 0.
 In \cite{BER97}, it is shown that there are only
two  automorphisms of $M$ at 0, namely $H_1(z,w)=(z,w)$ and $H_2(z,w)=(-z,w)$.
Thus, $ \C \big\{
j_0^K\big(\text{Aut}(M,0)\big)\big\}$ also has only two
elements. However, the group of holomorphic Segre preserving automorphisms of $\mathcal{M}$ at 0
(which according to Corollary \ref{thm:cor1.3} necessarily consists of a finite number of elements)
contains at least four maps: $$\HH_1(z,w,\chi,\tau)=(z,w,\chi,\tau) ,$$
$$\HH_2(z,w,\chi,\tau)=(-z,w,-\chi,\tau),$$ $$\HH_3(z,w,\chi,\tau)=(-z,-w,\chi,-\tau),$$
\begin{equation*}
\HH_4(z,w,\chi,\tau)=(z,-w,-\chi,-\tau).
\end{equation*}
\end{ex4.2}

In the next two examples, we will compare $\eta_0^K \big(\Aut\big)$ and $\C \big\{j_0^K\big(\text{Aut}(M,0)\big)\big\}$
for the family $\mathcal{F}$ given by
\begin{equation*}
\mathcal{F} = \Big\{M = \big\{\text{Im }w=c_1|z|^{2m}+c_2|z|^{2n} \big\} \, \Big| \,  1<m < n,  |c_1|^2+|c_2|^2 \neq 0 \Big\}  .
\end{equation*}
(We  exclude the \emph{Levi flat} case, $M=\{\text{Im }w=0\}$, as there is no finite jet determination for this $M$.)
Notice that each submanifold in $\mathcal{F}$ is of finite type and finitely \emph{degenerate} at 0.
%For $K = 1$, let $(\Lambda_1,\Lambda_2,\Lambda_3,\Lambda_4$,
%$\widetilde{\Lambda}_1,\widetilde{\Lambda}_2,\widetilde{\Lambda}_3, \widetilde{\Lambda}_4)$
% be  coordinates on $G_0^1(\C^2) \times G_0^1(\C^2)$, where
% $\Lambda_1$ (resp., $\widetilde{\Lambda}_1$) corresponds to $f_z(0)$ (resp., $\f_\X(0)$),
%$\Lambda_2$ (resp., $\widetilde{\Lambda}_2$)  corresponds to $g_w(0)$ (resp., $\g_\T(0)$),
%$\Lambda_3$ (resp., $\widetilde{\Lambda}_3$)  corresponds to $f_w(0)$ (resp., $\f_\T(0)$), and
%$\Lambda_4$ (resp., $\widetilde{\Lambda}_4$)  corresponds to $g_z(0)$ (resp., $\g_\X(0)$). Similarly, for $K = 2$, we will
%let $(\Lambda_1,\Lambda_2,\Gamma_3,\Gamma_4, \Lambda'',\widetilde{\Lambda}_1,\widetilde{\Lambda}_2,\widetilde{\Gamma}_3,
%\widetilde{\Gamma}_4, \widetilde{\Lambda}'')$
%be  coordinates on $G_0^2(\C^2) \times G_0^2(\C^2)$, where $\Lambda_1$, $\Lambda_2$, $\widetilde{\Lambda}_1$, and $\widetilde{\Lambda}_2$
%are as defined above,
% $\Gamma_3$ (resp., $\widetilde{\Gamma}_3$) corresponds to $f_{zw}(0)$ (resp., $\f_{\X\T}(0)$),
% $\Gamma_4$ (resp., $\widetilde{\Gamma}_4$) corresponds to $g_{ww}(0)$ (resp., $\g_{\T\T}(0)$),
% and $\Lambda''$ and $\widetilde{\Lambda}''$ represent the remaining coordinates.

\begin{ex4.3}   \label{ex4.3bb}
Assume   $c_1 \neq 0$ and $c_2=0$.  It can be shown (see \cite{Angle} for the calculations)
that any holomorphic Segre preserving automorphism of $\M$ at 0 is given by
\begin{equation*}
\mathcal{H}(z,w,\X,\T) =
\big(f(z,w),g(z,w),\f(\X,\T),\g(\X,\T)\big)
\end{equation*}
\begin{equation}  \label{eq:fg2}
 =  \left( \frac{az}{\sqrt[m]{1 + \alpha w}} \,\, , \,\, \frac{a^m\tilde{a}^m w}{1+\alpha w} \,\, , \,\,
\frac{\tilde{a}\X}{\sqrt[m]{1 + \alpha \T}} \,\, , \,\, \frac{a^m\tilde{a}^m \T}{1+\alpha \T} \right) ,
\end{equation}
where $a,\tilde{a} \in \C \backslash \{0\}$, $\alpha \in \C$, and $f$ and $\tilde{f}$ are expressed in terms of any branch of the $m^{\text{th}}$ root.

It immediately follows that any automorphism of $M$ at 0 is of the form
\begin{equation}  \label{eee}
H(z,w) =
%\big(f(z,w),g(z,w)\big) =
\left( \frac{az}{\sqrt[m]{1 + \alpha w}} \, , \, \frac{|a|^{2m} w}{1+\alpha w} \right),
\end{equation}
where $a \in \C \backslash \{0\}$, $\alpha \in \mathbb{R}$, and $f$ is expressed in terms of any branch of the $m^{\text{th}}$ root.

In this case, $\eta_0^2 \big(\Aut\big) = \C \big\{j_0^2\big(\text{Aut}(M,0)\big)\big\}$. Indeed, from (\ref{eee}), we see
that $\C \big\{j_0^2\big(\text{Aut}(M,0)\big)\big\}$ is given by
\begin{equation*}
\big\{\Lambda^g_w=\widetilde{\Lambda}^{\g}_\T=(\Lambda^f_z)^m (\widetilde{\Lambda}^{\f}_\X)^m \, , \,
\Lambda^g_{ww} = 2m (\Lambda^f_z)^{m-1}(\widetilde{\Lambda}^{\f}_\X)^{m} \Lambda^f_{zw} \, , \,
\widetilde{\Lambda}^{\g}_{\T\T} = 2m (\widetilde{\Lambda}^{\f}_\X)^{m-1}(\Lambda^f_z)^m \widetilde{\Lambda}^{\f}_{\X\T} \, , \,
\end{equation*}
\begin{equation} \label{rrrrrr}
\Lambda^f_w = \Lambda^f_{ww}=\Lambda^f_{zz}=\Lambda^g_z = \Lambda^g_{zw}=\Lambda^g_{zz}
=\Lambda^{\f}_\T = \Lambda^{\f}_{\T\T}=\Lambda^{\f}_{\X\X}=\Lambda^{\g}_\X = \Lambda^{\g}_{\X\T}=\Lambda^{\g}_{\X\X}      =0 \big\} .
\end{equation}
It follows from (\ref{eq:fg2})  that $\eta_0^2 \big(\Aut\big)$ is also given by (\ref{rrrrrr}).
\end{ex4.3}

\begin{ex4.4}   \label{ex4.4bb}
Assume $c_1,c_2 \neq 0$.
It can be shown (see \cite{Angle} for the calculations) that
any holomorphic Segre preserving automorphism of $\M$ at 0
is of one of the following $n-m$ forms:
\begin{equation}    \label{e1bb}
\HH_c(z,w,\X,\T) = \left( a z, c^m w,  \frac{c}{a} \X , c^m \T \right) ,
\end{equation}
where $a \in \C \backslash \{0\}$ and $\displaystyle{  c \in \left\{ e^{ \frac{2i\pi r}{n-m} }: r=0, \ldots, n-m-1 \right\} }$
(i.e., $c$ is a primitive $(n-m)^{\text{th}}$ root of unity).

It immediately follows that any automorphism of $M$ at 0 is of the form:
\begin{equation}   \label{lasteq}
H(z,w) = (e^{i\theta} z, w) ,
\end{equation}
where $\theta \in \mathbb{R}$.

Thus, we see from (\ref{lasteq}) that $\C \big\{j_0^1\big(\text{Aut}(M,0)\big)\big\}$ is given by
\begin{equation*}
\big\{ \Lambda^g_w=\widetilde{\Lambda}^{\g}_\T=1 \, , \, \Lambda^f_z\widetilde{\Lambda}^{\f}_\X = 1 \, , \, \Lambda^f_w=\widetilde{\Lambda}^{\f}_\T
= \Lambda^g_z=\widetilde{\Lambda}^{\g}_\X=0 \}
\end{equation*}
and thus has positive dimension.
For $n=m+1$, (\ref{e1bb}) implies that $\C \big\{j_0^1\big(\text{Aut}(M,0)\big)\big\} = \eta_0^1 \big(\Aut\big)$.
For $n>m+1$, however,  $\C \big\{j_0^1\big(\text{Aut}(M,0)\big)\big\} \varsubsetneqq \eta_0^1 \big(\Aut\big)$.
Indeed, we see from (\ref{e1bb}) that $\eta_0^1 \big(\Aut\big)$ is equal to the disjoint union of exactly
$n-m$ distinct cosets of $\C \big\{j_0^1\big(\text{Aut}(M,0)\big)\big\}$.
\end{ex4.4}

\end{document}